\documentclass[11pt]{ascarticle}

\include{preamble}

\newcommand{\suppdoi}{https://dx.doi.org/10.5281/zenodo.7128677}

\newcommand{\fmodesvg}{0}
\newcommand{\fmodepdf}{1}
\newcommand{\figuremode}{\fmodepdf}

\DeclareUnicodeCharacter{2212}{-}

\newcommand{\sitext}[1]{\shortintertext{#1}}
\newcommand{\siand}{\shortintertext{and}}

\DeclareMathOperator\Div{Div}
\AtBeginDocument{ % for compatiblity with XeLaTeX's unicode-math
  \let\div\relax
  \DeclareMathOperator{\div}{div}
}

\DeclareMathOperator\Curl{Curl}
\DeclareMathOperator\curl{curl}

\DeclareMathOperator\Grad{Grad}
\DeclareMathOperator\grad{grad}

\DeclareMathOperator\Tr{Tr}
\DeclareMathOperator\Det{Det}

\newcommand{\solf}[1]{\hat{#1}}

\newcommand{\inff}[1]{\underset{#1}{\vphantom{sup}\inf}}
\newcommand{\supp}[1]{\underset{#1}{\vphantom{sup}\sup}}

\newcommand{\pdiff}[2]{\frac{\partial#1}{\partial#2}}

\newcommand{\norm}[1]{\lVert #1 \rVert}

\newcommand{\mcolon}{:}

\newcommand{\trate}[1]{\kern2.5pt{\dot{\kern-2.5pt{#1}}}}
\newcommand{\pert}[1]{\kern2.2pt{\widetilde{\kern-2.2pt{#1}}}}
\newcommand{\mic}[1]{\pert{#1}}

\newcommand{\mbf}[1]{\mathbf{#1}}

\newcommand{\bzero}{{\mbf{0}}}
\newcommand{\bone}{{\mbf{1}}}

\newcommand{\bA}{{\mbf{A}}}
\newcommand{\bB}{{\mbf{B}}}
\newcommand{\bC}{{\mbf{C}}}

\newcommand{\bF}{{\mbf{F}}}
\newcommand{\bFT}{{\mbf{F}^\text{T}}}
\newcommand{\bFi}{{\mbf{F}^{-1}}}
\newcommand{\bFiT}{{\mbf{F}^{-\text{T}}}}

\newcommand{\bH}{{\mbf{H}}}

\newcommand{\bN}{{\mbf{N}}}
\newcommand{\bP}{{\mbf{P}}}

\newcommand{\ba}{{\mbf{a}}}
\newcommand{\bb}{{\mbf{b}}}

\newcommand{\bbf}{{\mbf{f}}}

\newcommand{\bh}{{\mbf{h}}}

\newcommand{\bbm}{{\mbf{m}}}
\newcommand{\bn}{{\mbf{n}}}

\newcommand{\bt}{{\mbf{t}}}

\newcommand{\Bsigma}{{\boldsymbol{\upsigma}}}
\newcommand{\Bsigmamw}{{\Bsigma^{\text{MW}}}}

\newcommand{\rd}{\mathop{}\!\mathrm{d}}

\newcommand{\dv}{\rd v}
\newcommand{\dV}{\rd V}
\newcommand{\da}{\rd a}
\newcommand{\dA}{\rd A}

\newcommand{\perunit}[1]{\,[\si{#1}]}
\newcommand{\perarbitraryunit}{\,[\textit{a.u.}]}

\newcommand{\unitB}{\tesla}

\DeclareSIUnit{\unitone}{1}

\newcommand{\Psivac}{\Psi^\text{vac}}

\newcommand{\psib}{\Psi_t^b}
\newcommand{\psih}{\Psi_t^h}
\newcommand{\psibvac}{\Psi_t^{\text{vac},b}}
\newcommand{\psihvac}{\Psi_t^{\text{vac},h}}
\newcommand{\psibmat}{\Psi_t^{b}}
\newcommand{\psihmat}{\Psi_t^{h}}
\newcommand{\Psib}{\Psi_0^b}
\newcommand{\Psih}{\Psi_0^h}
\newcommand{\PsiCB}{\Psi_0^B}
\newcommand{\PsiCH}{\Psi_0^H}
\newcommand{\PsiCBvac}{\Psi^{\text{vac},B}_0}
\newcommand{\PsiCHvac}{\Psi^{\text{vac},H}_0}

\newcommand{\bodyforce}{\boldsymbol{f}^\text{b}}
\newcommand{\bfmagb}{\boldsymbol{f}^{\text{mag},b}}
\newcommand{\bfmagh}{\boldsymbol{f}^{\text{mag},h}}
\newcommand{\PKone}{\bP}
\newcommand{\PKonemw}{{\PKone^{\text{MW}}}}
\newcommand{\Bsigmamagb}{\Bsigma^{\text{mag},b}}
\newcommand{\Bsigmamagh}{\Bsigma^{\text{mag},h}}

\newcommand{\px}{x}
\newcommand{\pX}{X}

\newcommand{\baemic}{\mic{\ba}\kern0.5pt{^\text{e}}\kern-0.5pt}

\newcommand{\bammic}{\mic{\ba}\kern0.5pt{^\text{m}}\kern-0.5pt}

\newcommand{\Fhat}{\kern2.5pt{\widehat{\kern-2.5pt\bF}}}
\newcommand{\msat}{m^\text{s}}
\newcommand{\defomap}{\boldsymbol{\varphi}}

\newcommand{\vol}{\mathcal{B}}

\newcommand{\volb}{\vol}

\newcommand{\volf}{\vol^\infty}
\newcommand{\volfree}{\vol^\text{a}}

\newcommand{\boun}{\partial\vol}
\newcommand{\bounb}{\partial\volb}
\newcommand{\bounf}{\partial\volf}
\newcommand{\bounfree}{\partial\volfree}

\newcommand{\labelu}{displacement $u_2(0,R) / R$}

\newcommand{\labelmwnb}{\hspace*{-2mm}MW (V)}

\newcommand{\labeltcnb}{\hspace*{-2mm}TC (V)}

\newcommand{\labelnv}{\hspace*{-2mm}NV}

\newcommand{\SIshort}[2]{\SI[output-exponent-marker =
    \text{e},zero-decimal-to-integer]{#1}{#2}}

\begin{document}
    
    \unitlength1cm
    
    \begin{frontmatter}
        
    \title{Curing spurious magneto-mechanical coupling in soft non-magnetic materials}
    
    \author{Matthias Rambausek\corref{cor1}}
    \ead{matthias.rambausek@tuwien.ac.at}
    
    \author{Joachim Schöberl}
    \ead{joachim.schoeberl@tuwien.ac.at}
        
    \address{%
        Institute of Analysis and Scientific Computing, TU Wien, Wiedner Hauptstrasse 8-10, 1040 Wien, Austria}
    
    \cortext[cor1]{Corresponding author}

    \begin{abstract}
        The present work is concerned with the issue of spurious coupling 
        effects that are pervasive in fully coupled magneto-mechanical finite
        element simulations involving very soft non-magnetic or air-like media. 
        We first address the characterization of the spurious
        magneto-mechanical effects and their intuitive interpretation based on
        energy considerations.
        Then, as main contribution, we propose two new ways to prune the
        undesired spurious magneto-mechanical coupling in non-magnetic media.
        The proposed methods are compared with
        established methods in the context of magnetic bodies embedded in (i)
        air or vacuum and (ii)
        very soft elastic non-magnetic media.
        The comparison shows that the proposed approaches are accurate and effective.
        They, furthermore, allow for a consistent linearization of the coupled
        boundary value problems, which is crucial for the simulation of compliant
        structures. 
        For reproducibility and accessibility of the proposed methods, we provide our implementations with Netgen/NGSolve as well
        as all codes necessary for the reproduction of our results as supplementary
        material.
    \end{abstract}

\end{frontmatter}

\section{Introduction}
\label{sec:introduction}

In this work we address the fundamental numerical issue of spurious magnetic
coupling in full-field simulations of extremely compliant magnetic bodies and
structures. The bodies that we have in mind are composites consisting of an elastomer
matrix with almost rigid ferromagnetic inclusions, i.e., magnetorheological
elastomers (MREs). 
For experimental data on MREs we refer to 
\citet{jolly+etal1996},
\citet{bednarek1999},
\citet{ginder+etal2002},
\citet{martin+etal2006},
\citet{boese+roeder2009},
\citet{danas+kankanala+triantafyllidis2012},
\citet{bodelot+etal2018} and
\citet{garcia-gonzalez+etal2021}
as well as to the extensive review by
\citet{bastola+hossain2020}.
These findings have inspired many interesting applications of MREs such magnetically
actuated valves
\citep{boese+rabindranath+ehrlich2012}, magnetically tunable dampers
\citep{stewart+etal1998} and vibration isolators
\citep{ginder+etal2001,hitchcock+gordaninejad+fuchs2005,hitchcock+gordaninejad+fuchs2006} 
as well as micropumps \citep{cesmeci+etal2022}.
With matrix materials becoming softer and bodies and structures becoming more
compliant, new opportunities for
applications open up in medical and biological contexts
\citet{gonzalez-rico+etal2021,moreno-mateos+etal2022b}. However,
well-established experimental, numerical and theoretical approaches face new
challenges. For experimental data in this regard we refer to
\citet{nikitin+etal2004}, \citet{stepanov+etal2007,stepanov+etal2008}
and rather recently \citet{moreno+etal2021,moreno-mateos+etal2022c}, 
which set several opportunities for future research.

In this context, full-field simulations context are paramount for the understanding and
the prediction of the strongly coupled nonlinear magneto-mechanical response of
the materials and bodies of interest. The theoretical groundwork has been
established in the second half of the 20\textsuperscript{th} century
\citep{tiersten1964,brown1966,pao+hutter1975} and modernized roughly two
decades ago
\citep{kovetz2000, dorfmann+ogden2003,dorfmann+ogden2004a,dorfmann+brigadnov2004,kankanala+triantafyllidis2004,ericksen2006}.
With the availability of corresponding finite element methods%
\footnote{Note that the fundamental equations of finite electro-elasticity are
    formally very similar to those for finite magneto-elasticity. As a
    consequence, there was and is a fruitful exchange of computational methods
    between both fields.}%
, e.g.,
\citep{semenov+etal2006,vu+steinmann+possart2007,bustamante+dorfmann+ogden2011},
numerical simulations have been widely applied to augment purely theoretical as
well as experimental data and insight.
In this regard it is important to account for the composite nature of
magnetorheological elastomers with computational means as done by
\citet{javili+chatzigeorgiou+steinmann2013},
\citet{kalina+metsch+kaestner2016,kalina+etal2020},
\citet{miehe+vallicotti+teichtmeister2016},
\citet{keip+rambausek2016,keip+rambausek2017}, 
\citet{danas2017},
\citet{keip+sridhar2019}
and
\citet{kalina+etal2020}.
More theoretical or analytical approaches in this direction have been pursued,
e.g., by 
\citet{ponte-castaneda+galipeau2011},
\citet{galipeau+ponte-castaneda2013},
\citet{galipeau+ponte-castaneda2013b},
\citet{romeis+etal2017},
\citet{lefevre+etal2017},
\citet{romeis+toshchevikov+saphiannikova2019}
and
\cite{lefevre+etal2020}.
For explicit models with the capability to fit computational homogenization
results we mention \citet{mukherjee+etal2020} and \citet{kalina+etal2020}.
Models fitted against experimental data have been reported by 
\citet{danas+kankanala+triantafyllidis2012} and \citet{haldar2021}, 
whereby the latter takes viscoelastic effects into account.

A particular breed of MREs based on remanently magnetized particles dates back at
least to the exploratory experiments of \citet{stepanov+chertovich+kramarenko2012}.
Such ``hard'' MREs (h-MREs) have recently received increased
thanks to advances in additive manufacturing of h-MRE \emph{structures}
\citep{kim+etal2018}. One prominent field of application of this new technology
is soft robotics \citep{kim+etal2019}.
Numerical modeling of h-MREs under assumption of a \emph{prescribed} remanent
magnetization has been addressed by
\citet{garcia-gonzalez2019},
\citet{zhao+etal2019},
\citet{garcia-gonzalez+hossain2021},
\citet{kadapa+hossain2022},
\citet{dadgar-rad+hossain2022}. The resulting models reduce to finite elasticity
with particular loading conditions.
In contrast to that, the important aspect of ferromagnetic evolution in the
particles, e.g. during the magnetization process of h-MREs, has been accounted
in the methodological contributions of \citet{kalina+etal2017} and
\citet{rambausek+mukherjee+danas2022}. Macroscopic material models that are
able to reproduce microstructural simulations of h-MREs have been developed by
\citet{mukherjee+rambausek+danas2021} and \citet{mukherjee+danas2022}.
For comprehensive reviews on h-MREs we refer to
\citet{lucarini+hossain+garcia-gonzalez2022} but also to
\citet{wu+etal2020} where both (magnetically) ``soft'' MREs
and h-MREs were considered.

\begin{figure}[!ht]
    \centering
    \ifthenelse{\equal{\figuremode}{\fmodesvg}}{%
        \includesvg[width=0.5\textwidth,pretex=%
        \newcommand{\lmagn}{magnetized}%
        \newcommand{\lnonmagn}{non-magnetic}%
        ]{from_thesis_ellipse_in_carrier_distortion_examples_a}}{%
        \includegraphics{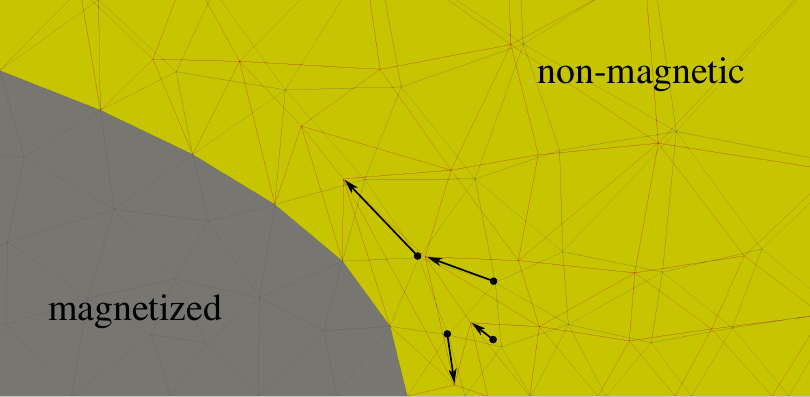}}%
    \caption{A spatially fixed, rigid and magnetized body (grey)
        embedded in a soft non-magnetic material (yellow).
        The rigid, magnetized body does not deform. However,
        \emph{spurious} deformation (red vs. grey mesh) is observed in the
        non-magnetic
        domain in response to magnetic fields 
        \citep[Chap.~9.3.1;Figure~9.17a]{rambausek2020}.}
    \label{fig:demo_of_spurious_defo}
\end{figure}

No matter whether ``soft'' or ``hard'' MREs are considered in simulations, one
almost always encounters the pervasive issue that will be the addressed by
this contribution: spurious magneto-mechanical coupling in sufficiently soft or
compliant magnetoactive bodies as depicted in
Figure~\ref{fig:demo_of_spurious_defo} \citep{rambausek2020}.
There one can observe excessive deformation in the non-magnetic material
(yellow) embedding a clamped, \emph{rigid} magnetized body (grey) where no
deformation at all is expected according to the laws of physics.
These effects are ubiquitous but only become relevant in
extremely soft materials. This can be for example the case when 
a sufficiently soft non-magnetic medium is employed to model the air
surrounding some magnetic specimen \citep{keip+rambausek2016}.
However, it may also be observed
when increasingly soft%
\footnote{Shear modulus in the range of \SI{10}{\kilo\pascal} or lower.}
matrix materials are employed as is the case for example
in the experiments of \citep{moreno+etal2021}.
While both situations have the commonality of an extremely soft non-magnetic
material being present in the boundary value problem, in the first case this is
purely for convenient modeling of a medium with no significant stiffness at
all, in the second case the soft material is physical reality. 
Due to that,
in the former case, one has more modeling and methodological freedom than in the
latter. In fact, \citet{pelteret+etal2016} proposed an elegant staggered
approach to account for air or vacuum surrounding a magnetic body that does not
make use of an auxiliary soft material and hence does not suffer from spurious
coupling in the present sense. However, their approach is not applicable
when a soft elastic material is indeed physically present. Also, due to its
partially decoupled nature and the inherent inconsistent linearization of the
discretized boundary value problem, 
it suffers from convergence issues in the case of highly compliant structures.
As a second method that successfully circumvents spurious coupling effects we
mention the staggered approach of \citet{liu+etal2020}. Their approach is based
on surface currents computed in a purely magnetic boundary value problem. 
In a second step, these currents are employed to compute magneto-mechanical
tractions on boundaries of magnetic bodies for the mechanical subproblem. 
Due to the decoupled magnetostatic and mechanical subproblems, the scheme
requires a certain number of interations between both of them to account for
the strong coupling and attain convergence. While the procedure is effective for the
analysis of slender magneto-active structures, its staggered nature and the
explicit use of electric currents for force computation render the scheme less
suited for adoption in existing fully coupled (monolithic) methods.
A third approach for the modeling of surrounding air in finite element
simulation are non-local constraints that bind the motion in the air domain to
the motion of the boundary of the magnetic bodies under consideration
\citep{psarra+etal2019,mukherjee+rambausek+danas2021}. While
this method does not introduce artificial stiffness and linearizes properly it
leads to a much denser linear system. Furthermore, it is difficult to implement
for general geometries and arrangements of multiple magnetic bodies.
For some additional insights on existing monolithic and staggered schemes
we refer to \citet{rambausek+mukherjee+danas2022}.

The methods proposed in this contribution exploit the non-magnetic nature of
the soft material (surrounding air or soft carrier/matrix material) to obtain
formulations that integrate well with existing magneto-elasticity frameworks,
i.e. it is not much more difficult to implement than the naive monolithic
schemes, for both the surrounding air and the soft matrix case. Furthermore,
the schemes obtained linearize properly and at the same
time do not suffer from spurious coupling and are thus promising candidates for
widespread adoption. To facilitate this process we provide our implementation
with Netgen/NGSolve \citep{schoeberl1997,schoeberl2014,ngsolveweb} as well as codes
for our numerical examples as \href{\suppdoi}{supplementary material}.

The remainder of this contribution is structured as follows.
In Section~\ref{sec:theoretical_background} we outline the fundamental equations
of finite magneto-elasticity as well as important considerations on
magneto-mechanical interactions. Section~\ref{sec:spurious_forces} then
presents the spurious magneto-mechanical coupling that is the motivation for
the present work and discusses the underlying general pattern. 
Then, in Section~\ref{sec:remedies} we arrive at the heart of our contribution where
we propose two approaches that eliminate
spurious magneto-mechanical interactions in non-magnetic domains.
We sketch numerical implementations of the schemes and assess their 
accuracy by a comparison with existing methods as well as their convergence under
mesh refinement. After that we showcase both proposed
approaches in a challenging example involving a magnetic solid and a very soft
non-magnetic solid under gravitational and magnetic loads embedded in an air-like
medium.
We close our contribution with the conclusion in Section~\ref{sec:conclusion}.

\section{Theoretical background}
\label{sec:theoretical_background}

\begin{figure}[!ht]
    \centering
    \ifthenelse{\equal{\figuremode}{\fmodesvg}}{%
    \includesvg[width=0.5\textwidth,pretex=%
    \small%
    \newcommand{\labelbinfty}{$\bb^\infty$}%
    \newcommand{\labeldomain}{$\volfree$; $\volf = \volfree \cup \left(\bigcup\volb^i\right)$}%
    \newcommand{\labelbc}{mech. support}%
    ]{bodies_in_free_space_simple}}{%
    \includegraphics{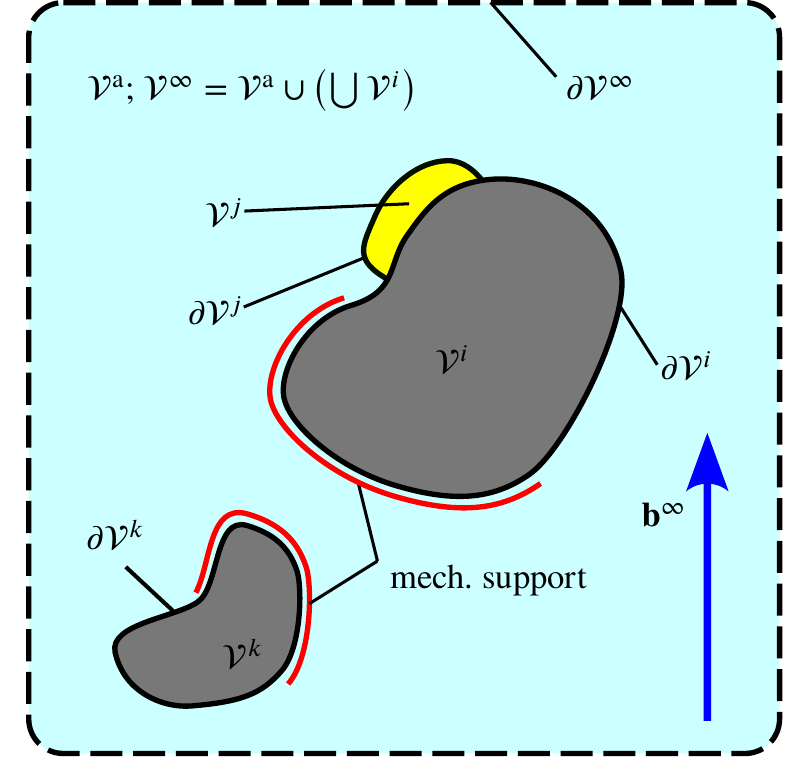}}%
    \caption{The overall setting considered in this work. Magnetic bodies (grey;
        $\volb^{i,k}$) and non-magnetic bodies (yellow; $\volb^{j}$) in ``free''
        space (light cyan; vacuum,
        air-like medium) exposed to a uniform external magnetic field
        $\bb^\infty$. We employ the symbol $\volf$ for the full domain including
        all sub-bodies and $\volfree$ for the ``remaining'' empty or air-like
        domain.}
    \label{fig:bodies_in_free_space}
\end{figure}

For the problem of spurious magneto-mechanical coupling investigated in this
work, we restrict ourselves to quasi-static finite elasticity 
coupled to magnetostatics in absence of free currents.
We first recall the fundamental theory for the mechanical part, 
then the magnetostatic and finally bring both fields together.
Figure~\ref{fig:bodies_in_free_space} provides a overview of the setting
considered. 
There we depict magnetic bodies (grey; $\volb^{i,k}$) and non-magnetic bodies
(yellow; $\volb^{j}$) embedded in ``empty'' space $\volfree$ (light cyan; vacuum,
air-like medium).
The bodies are exposed to a uniform external magnetic field $\bb^\infty$ and may
have mechanical support. While not indicated, we may of course consider mechanical
body forces like gravity as well. The ``full'' domain consisting of $\volf$ and
all (solid) bodies $\volb^i$ is denoted as $\volf$. 
The (computational) outer boundary $\bounf$ is indicated by a dashed line. 
Physically, there is
no boundary but in finite element simulations the
infinite space domain is truncated%
\footnote{An alternative to truncation of the computational domain is to couple
    finite element
    simulations of the solid bodies
    with boundary elements employed for the solution of the magnetic field
    equations in the ``empty'' domain \citep{vu+steinmann2010,vu+steinmann2012}.} 
at a larger distance from the bodies under
consideration
\citep{bustamante+dorfmann+ogden2011,miehe+ethiraj2012,keip+rambausek2016,schroeder+etal2022}.
For definiteness, we mechanically fix the outer boundary $\bounf$.

\subsection{Quasi-static finite hyperelasticity}
\label{sec:into_elasticity}

In the context of finite strains we have to distinguish between the
\emph{Eulerian} or \emph{current configuration} and the \emph{Lagrangian} or
\emph{reference\footnote{In the present context, one may regard the
         \emph{initial configuration}, i.e. undeformed bodies and media as at
         the beginning of the deforming process, as
         reference.}
configuration}.
The former simply refers to the (deformation) state of media as an
external observed perceives them, whereas the latter takes the perspective of
material points. Following a common convention, we denote quantities in the
Eulerian configuration by lowercase symbols, whereas uppercase symbols are
employed for the Lagrangian configuration.
In some instances, we employ subscript ``0'' for Lagrangian and subscript ``t''
for Eulerian quantities for definiteness.
The Eulerian and Lagrangian configurations are connected through the
map $\defomap$, that assigns to each Lagrangian position $\pX$ (a material
point) a Eulerian position $\px = \defomap(\pX)$. 
The corresponding tangent map or Jacobian $\bF = \pdiff{\defomap}{\pX}$. 
In Cartesian coordinates $\bF$ is often written as 
\begin{equation}
\label{eq:defograd}
\bF(\pX) = \Grad\defomap(\pX) \quad\Rightarrow\quad \Curl\bF\equiv0
\end{equation} and
commonly referred to as deformation gradient.
Its determinant $J=\Det\bF$ is a measure for the change of volume elements
\begin{align}
    \label{eq:change_of_volume}
    \dv = J\dV,
\end{align}
where $\dv$ is the Eulerian and $\dV$ the Lagrangian volume element.

On the kinetic side, the balance of momentum in a quasi-static setting reads
\begin{align}
\label{eq:balance_of_momentum}
\div\Bsigma + \bodyforce &= 0 \quad\text{in}\: \volb_t
&\Leftrightarrow && \Div \PKone + \bodyforce_0 &= 0 
\quad\text{in}\: \volb_0
\end{align}
in the Eulerian and Lagrangian setting, respectively, where $\Bsigma$ denotes
the Cauchy stress, $\PKone$ the ``first Piola-Kirchhoff stress'' and
$\bodyforce$ denotes the body force density per Eulerian volume.
The surface-integral-preserving Piola-type transform relating $\Bsigma$ and $\PKone$
is 
\begin{equation}
    \label{eq:transform_of_stress}
    \Bsigma(\pX) = J^{-1} \PKone \cdot \bFT.
\end{equation}
The Eulerian body force density $\bodyforce$ transforms as
\begin{equation}
    \label{eq:transform_of_body_force}
    \bodyforce(\pX) = J^{-1} \bodyforce_0,
\end{equation}
The corresponding boundary conditions are given as 
\begin{align}
\label{eq:traction_boundary_conditions}
\Bsigma\cdot\bn &= \bt \quad\text{on}\: \bounb_t
&\Leftrightarrow && \PKone \cdot \bN &= \left(J\,\sqrt{(\bFT \cdot \bN) \cdot (\bFT \cdot \bN)}\right) \bt \quad\text{on}\: \bounb_0,
\end{align}
where $\bt$ is a given Eulerian traction, i.e. a surface force density per
Eulerian area and $\bn$ is the outward-pointing unit normal. 
The complicated factor on the right is the ratio of area elements
$\da/\dA$ which transforms the Eulerian to a Lagrangian traction.

In the context of hyperelasticity, the stresses are obtained as partial
derivatives of a strain energy density per mass $\psi(\bF)$
\begin{align}
    \label{eq:constitutive_relation_hyperelasticity}
    &\Bsigma = \rho \pdiff{\psi}{\bF}
    & \text{and} &
    &\PKone = \rho_0  \pdiff{\psi}{\bF} \cdot \bFT
\end{align}
where $\rho$ and $\rho_0$ are the mass densities per current (Eulerian) and
initial (Lagrangian) volume. The latter are related through
\begin{align}
    \label{eq:transform_of_mass_density}
    \rho = J^{-1}\rho_0 .
\end{align}

As the final modeling ingredient we introduce the
strain energy density per Lagrangian $\Psi_0(\bF) = \rho_0\psi(\bF)$
volume, where $\rho_0$ is the respective mass density.
Correspondingly, the strain energy density per Eulerian volume is
$\Psi_t(\bF) = \rho\psi(\bF)$.

The balance of momentum is accompanied
the balance of moment of momentum expressed as
\begin{equation}
\label{eq:balance_of_moment_of_momentum}
\Bsigma^{\text{T}} = \Bsigma ,
\end{equation}
i.e. the symmetry of the \emph{Cauchy} stress.
In the present setting, \eqref{eq:balance_of_moment_of_momentum} is ensured by an
objective (material frame indifferent)
strain energy density. Material frame indifference is achieve, for example, 
by parameterizing the strain energy density in terms of the ``right Cauchy-Green''
tensor 
\begin{equation}
\label{eq:right_Cauchy_Green}
\bC = \bFT \cdot \bF
\end{equation}
such that $\Psi(\bF)_0 = \Psi_0(\bC(\bF))$ which yields
\begin{align}
\label{eq:Cauchy_symmetric}
\Bsigma = J^{-1} \pdiff{\Psi_0}{\bF} \cdot \bFT = 2 \, J^{-1}
\bF \cdot \pdiff{\Psi_0}{\bC} \cdot \bFT .
\end{align}
The last term clearly is symmetric and thus $\Bsigma$ is symmetric too.

Considering the case where the external volume and surface force densities are
\textit{conservative}, we introduce the respective loading potential 
densities (per Eulerian volume) $l^\bbf$ and $l^\bt$ via
\begin{align}
    \label{eq:mechanical_loading_potentials}
    &\bodyforce = -\pdiff{l^\bbf}{\defomap}
    & \text{and} &
    &\bt = -\pdiff{l^\bt}{\defomap} .
\end{align}
With these at hand, we may formulate a variational principle for finite hyperelasticity
\begin{align}
\label{eq:mechanical_energy_principle_Eulerian}
\hat\defomap = \arg \inff{\defomap}\left\{
\int_{\volb_t} \rho\psi(\bC) \dv
+ \int_{\volb_t} l^\bbf \dv 
+ \int_{\bounb_t} l^\bt \da
\right\},
\end{align}
where we remark that the integrals can be individually transformed to the 
Lagrangian configuration for computational convenience
with the help of \eqref{eq:change_of_volume},
\eqref{eq:transform_of_mass_density}
and the area transform factor in \eqref{eq:traction_boundary_conditions}.
For example, the fully transformed version reads
\begin{align}
\label{eq:mechanical_energy_principle_Lagrangian}
\hat\defomap = \arg \inff{\defomap}\left\{
\int_{\volb_0} \rho_0\psi(\bC) \dV
+ \int_{\volb_0} l_0^\bbf \dV
+ \int_{\bounb_0} l_0^\bt \dA
\right\} .
\end{align}

\subsection{Magnetostatics}
\label{sec:intro_magnetostatics}

The magnetostatic field equations reduced to the scope of the present work
read
\begin{subequations}
    \label{eq:magnetostatics}
\begin{align}
    \label{eq:magnetostatics_b}
    \div \bb &= 0 \quad\text{in}\: \volf_t & {\bb}\cdot\bn &= {\bb^\infty}\cdot\bn \quad\text{on}\:\bounf_t \\%& \bb &= \curl \ba \\
    \siand
    \label{eq:magnetostatics_h}
    \curl \bh &= \bzero \quad\text{in}\: \volf_t & {\bh} \times \bn &= {\bh^\infty} \times \bn \quad\text{on}\:{\bounf_t} ,%& \bh &= -\grad \phi,
\end{align}
\end{subequations}
where $\bb$ is the magnetic $b-$field (magnetic flux), $\bh$ is magnetic
$h-$field (magnetic intensity).
From \eqref{eq:magnetostatics} it is clear that $\bb$ and $\bh$ can be expressed
in terms of potentials as
\begin{align}
    \label{eq:magnetic_potentials}
    &\bb = \curl\ba & \text{and} && \bh =-\grad\phi ,
\end{align}
with $\ba$ being the (magnetic) vector potential and $\phi$ the auxiliary
\emph{scalar} magnetic potential.
While not of real physical significance, $\phi$ is particularly
convenient for computations.

It is possible to show that the equations introduced so far can be obtained
from two equivalent variational principles. The first, to which we refer as
magnetic energy principle reads
\begin{align}
    \label{eq:magnetic_energy_principle}
    \solf{\ba} = \arg \inff{\ba} \left\{ \int_{\volf} \psib(\px, \bb=\curl\ba) \dv - \int_{\bounf_t} (\bh^\infty\times\bn)\cdot\ba \da \right\}
\end{align}
where $\psib(\px,\bb)$ is the magnetic energy density per Eulerian volume, a
quantity to be specified depending on the medium occupying the respective point
$\px$ in space. At a solution $\hat{\ba}$ we have
\begin{align}
    \label{eq:magnetic_energy_principle_stat}
    \curl \bh &= 0 & \text{with} &&  \bh = \pdiff{\psib}{\bb} .
\end{align}

Dual to that we have the magnetic \emph{co}-energy principle
\begin{align}
    \label{eq:magnetic_co-energy_principle}
    \solf{\phi} = \arg \supp{\phi} \left\{ \int_{\volf_t} \psih(\px, \bh=-\grad\phi) \dv - \int_{\bounf_t} (\bb^\infty \cdot \bn)\,\phi\da \right\}
\end{align}
where $\psi(\bh)$ is the magnetic co-energy density, which is obtained from
$\psi(\bb)$ through a Legendre-Fenchel transform
\begin{align}
    \label{eq:legendre_transform}
    \psih(\bh) = \inff{\bb} \left\{-\bb \cdot \bh + \psib(\bb) \right\}.
\end{align}
At a solution $\hat{\phi}$ we have
\begin{align}
    \label{eq:magnetic_co-energy_principle_stat}
    \div \bb &= 0 & \text{with} &&  \bb = -\pdiff{\psih}{\bh} .
\end{align}

\subsection{Coupled magnetoelasticity}
\label{sec:intro_magnetoelasticity}

In order to couple magnetostatics
to finite elasticity we establish the following transforms
\begin{subequations}
        \label{eq:magneto-statics_Eulerian_Lagrangian}
\begin{align}
    &\phi(\pX) = \phi(\px) \circ \defomap(X)
    \\
    &\bh(\pX) = \grad(\phi(\defomap(\pX))) = \bFT \cdot \bH(\pX)
    & \text{with} &
    &\bH = \Grad(\phi(\pX))
    \\
    &\psih(\bh(\pX), \pX) = J^{-1} \Psih(\bh=\bFT\cdot\bH, \pX).
    \\
    &\ba(\pX) = \bFT \cdot \bA(\pX)
    \\
    &\bb(\pX) = \curl(\ba(\pX)) = J^{-1} \bF \cdot \bB(\pX)
    & \text{with} &
    &\bB = \Curl(\bA(\pX))
    \\
    &\psib(\bb(\pX), \pX) = J^{-1} \Psib(\bb = J^{-1}\bF\cdot\bB, \pX).
\end{align}
\end{subequations}
We note that the transforms involving $\bF^{-T}$ preserve line integrals,
those involving $J^{-1} \bF$ preserve surface integrals and those
involving $J^{-1}$ preserve volume integrals.

In what follows, we exploit that any material-frame indifferent $\Psih(\bh=\bFT\cdot\bH, \pX)$ can be recast to
the form $\PsiCH(\bC, \bH, \pX)$ and 
$\Psib(\bb = J^{-1}\bF\cdot\bB, \pX)$ to $\PsiCB(\bC, \bB, \pX)$
\citep{mukherjee+rambausek+danas2021,mukherjee+danas2022}.
The parameterizations of $\PsiCH$ and $\PsiCB$ have the advantage that they
directly generalize to the magneto-mechanical case since they depend
on both (mechanical) kinematic and magnetic quantities.
%
%With help of \eqref{eq:magneto-statics_Eulerian_Lagrangian} we reformulate
%the two variational principles \eqref{eq:magnetic_energy_principle} and  
%\eqref{eq:magnetic_co-energy_principle} to obtain
%\begin{align}
%    \label{eq:magnetic_energy_principle_Lagrangian}
%    &
%    \solf{\bA} = \arg \inff{\bA} \left\{ \int_{\mathbb{R}^3} \Psib(\pX, \bB=\Curl\bA) \dV \right\}
%    & \Rightarrow && \Curl\bH &= 0 && \text{with} && \bH = \pdiff{\Psib(\bB)}{\bB}
%    \\
%    \label{eq:magnetic_co-energy_principle_Lagrangian}
%    &\solf{\phi} = \arg \supp{\phi} \left\{ \int_{\mathbb{R}^3} \Psih(\pX, \bH=-\Grad\phi) \dV  \right\}
%    & \Rightarrow && \Div\bB &= 0 && \text{with} && \bB = -\pdiff{\Psih(\bH)}{\bH}
%\end{align}
%where we note the relation $\dv = J \dV$.
%
%We now generalize the abstract energy densities introduced above
%to a coupled energy density $\Psi(\bC, \bB)$ and a coupled
%``energy--co-energy'' density $\Psi(\bC, \bH)$, which may depend on both
%(mechanical) kinematic and magnetic quantities.
With these at hand, we combine the mechanical
\eqref{eq:mechanical_energy_principle_Lagrangian} and magnetic principles
\eqref{eq:magnetic_energy_principle} and
\eqref{eq:magnetic_co-energy_principle} to
the coupled variational principles \citep{bustamante+dorfmann+ogden2008}
\begin{align}
    \label{eq:magnetomechanical_energy_principle}
    \left\{\solf{\defomap}, \solf{\bA}\right\} = \arg \inff{\defomap,\bA}\left\{
        \int_{\volf_0}\PsiCB(\bC, \bB) \dV
        - \int_{\bounf_0} (\bh^\infty \times \bN)\cdot\bA\dA
        +\int_{\volb_0}l^\bbf_0 \dV
        +\int_{\boun_0}l^\bt_0 \dV
    \right\}
\siand
    \label{eq:magnetomechanical_energy_co-energy_principle}
    \left\{\solf{\defomap}, \solf{\phi}\right\} = \arg \inff{\defomap}\supp{\phi}\left\{
        \int_{\volf_0} \PsiCH(\bC, \bH) \dV
        - \int_{\bounf_0} (\bb^\infty \cdot \bN)\,\phi\dA
        +\int_{\volb_0}l^\bbf_0 \dV
        +\int_{\boun_0}l^\bt_0 \dV
    \right\} ,
\end{align}
whereby we note that, as a consequence
of fixing the outer boundary $\bounf$, the loading terms on that boundary in 
\eqref{eq:magnetomechanical_energy_principle} and
\eqref{eq:magnetomechanical_energy_co-energy_principle} indeed take this form with a
slight abuse of notation.
Moreover, outside $\volb$ the energy densities $\PsiCB(\bC, \bB)$
and $\PsiCB(\bC, \bH)$ reduce to the respective vacuum energy densities discussed in 
Subsection~\ref{sec:non-magnetic_media}.

\subsection{Specific relations for non-magnetic media}
\label{sec:non-magnetic_media}

In empty space (vacuum) and, to very good approximation, in all (technically)
\emph{non-magnetic}%
\footnote{Diamagnetism is commonly ignored in the wider context of this
    contribution.}
media the fields $\bh$ and $\bb$ are related through
%\footnote{The simplicity of this relation, in
%particular the scalar appearance of $\mu_0$, hides the fact that $\bh$ and
%$\bb$ are of different geometric nature such that more ingredients are required
%to render the equation \emph{tensorial} or \emph{covariant} as will be made
%clear below.}
\begin{equation}
    \label{eq:bh_relation}
    \bb = \mu_0\,\bh 
\end{equation}
with $\mu_0$ being the vacuum permeability.
The corresponding (co-)energy densities are given as
\begin{align}
    \label{eq:vacuum_energy_density}
    \psibvac(\bb) = \frac{1}{2\mu_0} \bb\cdot\bb 
    \siand
    \label{eq:vacuum_co-energy_density}
    \psihvac(\bh) = -\frac{\mu_0}{2} \bh\cdot\bh,
\end{align}
from which \eqref{eq:bh_relation} can be easily computed via
\eqref{eq:magnetic_energy_principle_stat}\textsubscript{2}
or
\eqref{eq:magnetic_co-energy_principle_stat}\textsubscript{2},
respectively.

An important consideration in this context is that non-magnetic media do not
experience any ``magnetic forces''. This means that they will not move or
deform under magnetic field. Also, in general they remain non-magnetic under
deformation.

This is reflected by the following observation: Consider any of the energy
densities \eqref{eq:vacuum_energy_density} or    
\eqref{eq:vacuum_co-energy_density} and employ them in the respective 
magnetomechanical variational principle
\eqref{eq:magnetomechanical_energy_principle} or
\eqref{eq:magnetomechanical_energy_co-energy_principle}.
For this purpose, the densities are expressed in their Lagrangian form
\begin{align}
    \label{eq:vacuum_energy_density_Lagrangian}
    \PsiCBvac(\bC, \bB) = \frac{1}{2 J \mu_0} \bB\cdot(\bC\cdot\bB)
    \siand
    \label{eq:vacuum_co-energy_density_Lagrangian}
    \PsiCHvac(\bC, \bH) = -\frac{J \mu_0}{2} \bH\cdot(\bC^{-1}\cdot\bH) .
\end{align}
Then, \emph{irrespective whether magnetostatic energy or co-energy is used}, the
resulting Cauchy-type stress field in this particular case is
\begin{align}
    \label{eq:cauchy-type_maxwell_stress_in_vacuum}
    \Bsigmamw := 2\,J^{-1}\bF\cdot\pdiff{\Psivac}{\bC}\cdot\bFT = \frac{1}{\mu_0} \bb \otimes\bb - \frac{1}{2\mu_0}\norm{\bb}\,\bone 
    = \bh \otimes\bb - \frac{1}{2} (\bh\cdot\bb) \bone 
    = {\mu_0} \bh \otimes\bh - \frac{\mu_0}{2} \norm{\bh}\,\bone,
\end{align}
which is commonly known as the
\emph{Maxwell stress}%
\footnote{We will employ ``Maxwell stress'' only for the particular forms below
    in vacuum. There is no unique notion of a Maxwell stress in matter such
    that we prefer avoid the resulting ambiguities.}.

The key properties of the Maxwell stress are 
(i) that its Eulerian form only
depends on magnetic quantities but no mechanical (kinematic) ones and 
(ii) that it is divergence-free for any magnetostatic solution in non-magnetic
media.
Combining both reveals that there cannot be any magnetic force densities in
such a medium irrespective of the deformation state.
An equivalent viewpoint is that Eulerian magnetostatic energy an co-energy are
invariant with respect to deformation in the interior of a non-magnetic medium.

\begin{rmk}
    \label{rmk:permeability_as_hodge_operator}
The simplicity of the vacuum constitutive relation \eqref{eq:bh_relation}, in
particular the scalar appearance of $\mu_0$, hides the fact that $\bh$ and
$\bb$ are of different geometric nature. By ``pulling-back'' this equation to
the Lagrangian configuration, we obtain
\begin{align}
    \label{eq:BH_relation}
    \bB = J \bFi \cdot (\mu_0 \bFiT) \cdot \bH = (J \bC^{-1} \mu_0) \cdot \bH,
\end{align}
which reveals the tensorial object
represented by $\mu_0$.
In terms of differential geometry, $\bh$ corresponds to a one-form whereas
$\bb$ corresponds to a two-form, which roughly means that the former is a
``gradient-type'' field and the latter is of ``curl-type''. Both can be
expressed as vectors, but they transform differently.
A tensorial object that relates fields of different character in the above
sense is an instance of a \emph{Hodge} operator.
\end{rmk}

\subsection{Magnetic media and magnetic forces}
\label{sec:magnetic_media}

In accordance with the above ``definition'' of a non-magnetic medium we regard
any medium for which \eqref{eq:bh_relation} does not hold in general as
\emph{magnetic}%
or
\emph{magnetized}.
The mismatch is connected to the magnetization $\bbm$ which appears as
\begin{align}
    \label{eq:magnetization}
    \bb = \mu_0 (\bh + \bbm),
\end{align}
such that non-magnetic media have $\bbm\equiv0$. 

It is important to note here that the generic
constitutive relations
\eqref{eq:magnetic_energy_principle_stat}\textsubscript{2} and
\eqref{eq:magnetic_co-energy_principle_stat}\textsubscript{2}
still hold. In particular \citep{bustamante+dorfmann+ogden2008},
\begin{align}
    \label{eq:constitutive_relation_for_magnetization_b}
    &\bbm = \pdiff{\psibvac(\bb)}{\bb} - \pdiff{\psib(\bF, \bb)}{\bb} 
          = \pdiff{\left(\psibvac(\bb) - \psib(\bF, \bb)\right)}{\bb}
          = -\pdiff{\psibmat(\bF, \bb)}{\bb}
\siand
    \label{eq:constitutive_relation_for_magnetization_h}
    &\mu_0\bbm = \pdiff{\psihvac(\bh)}{\bh} - \pdiff{\psih(\bF, \bh)}{\bh}
               = \pdiff{\left(\psihvac(\bh) - \psih(\bF, \bh)\right)}{\bh}
               = -\pdiff{\psihmat(\bF, \bh)}{\bh}, 
\end{align}
where $\psib(\bF, \bb)$ and $\psih(\bF, \bh)$ are \textit{total} energy
densities including the vacuum energy densities whereas
$\psibmat(\bF, \bb)$ and $\psihmat(\bF, \bh)$ are contributions from matter.

Concerning stresses and force densities due to the presence of magnetic fields,
we consider
\citep[Table~1]{kankanala+triantafyllidis2004}
\begin{align}
\label{eq:magnetic_stresses_b}
    &\Bsigmamagb     = \Bsigmamw + 
        J^{-1}
        \pdiff{\psibmat(\bF, \bb)}{\bb}\cdot\pdiff{J^{-1}\bF\cdot\bB}{\bF}
        \cdot
        \bFT 
        = \bh\otimes\bb - \frac{\mu_0}{2}(\norm{\bh} - \norm{\bbm})\bone
    \\
    \label{eq:magnetic_forces_b}
%    &
    &
    \bfmagb = -\div\Bsigmamagb = \bbm \cdot (\grad{\bb})
    \siand
\label{eq:magnetic_stresses_h}
    &\Bsigmamagh = \Bsigmamw + 
        J^{-1}
        \pdiff{\psibmat(\bF, \bh)}{\bh}\cdot\pdiff{\bFiT\cdot\bH}{\bF} 
        \cdot
        \bFT
        = \bh\otimes\bb - \frac{\mu_0}{2}\norm{\bh}\,\bone
    \\    
    \label{eq:magnetic_forces_h}
%    &
    &
    \bfmagh = -\div\Bsigmamagh = \mu_0(\grad{\bh})\cdot\bbm
        ,
\end{align}
where we note that results for the divergence were obtained for $\curl\bh=0$.
The fact that the stress contributions are obviously not divergence-free
can be interpreted as the presence of magnetic volume force
densities. The corresponding surface force densities at an interface to a
non-magnetic medium are \citep[Table~1]{kankanala+triantafyllidis2004}
\begin{align}
    \label{eq:magnetic_tractions_b}
    \left(\Bsigmamw\big|_{\text{out}} -
    \Bsigmamagb\big|_{\text{in}}\right)\cdot\bn
    = 
    \frac{-\mu_0}{2} \left[\left(\bbm\cdot\bn\right)^2 + \left(\bbm\cdot\bbm\right)\right] \bn
\siand
    \label{eq:magnetic_tractions_h}
    \left(\Bsigmamw\big|_{\text{out}} -
     \Bsigmamagh\big|_{\text{in}}\right)\cdot\bn
    = 
    \frac{-\mu_0}{2} \left(\bbm\cdot\bn\right)^2 \,\bn ,
\end{align}
where $\bn$ is the \emph{outward pointing normal vector}.

\begin{rmk}
    \label{rmk:symmetry_of_magnetic_stress}
    In magnetically anisotropic media, i.e. media where $\bb$ is not
    necessarily parallel to $\bh$, the magneto-mechanical
    Cauchy-type stresses \eqref{eq:magnetic_stresses_b} and
    \eqref{eq:magnetic_stresses_h} are not symmetric but only the 
    \emph{total} Cauchy-type stress \eqref{eq:Cauchy_symmetric} is.
    This corresponds to the presence of magnetic torques.
    Examples for such materials are MREs with anisotropic particle distributions
    \citep{boczkowska+awietjan2012,danas+kankanala+triantafyllidis2012}
    and h-MREs when remanently magnetized.
    We refer to
    \citet{mukherjee+rambausek+danas2021} for a related discussion of the
    latter case.
\end{rmk}

\begin{figure}[!ht]
    \centering
    \ifthenelse{\equal{\figuremode}{\fmodesvg}}{%
    \includesvg[width=0.5\textwidth,pretex=%
    \small%
    \newcommand{\lfields}{$\bb = \mu_0\bh$}%
    \newcommand{\lsurf}{$\mathcal{S}_t$}%
    \newcommand{\lmag}{$\bb = \mu_0(\bh + \bbm)$}%
    \newcommand{\ltraction}{$\Bsigmamw\cdot\bn$}%
    \newcommand{\lnormal}{$\bn$}%
    \newcommand{\lbody}{$\mathcal{B}_t$}%
    ]{body_in_free_space_surface}}{%
    \includegraphics{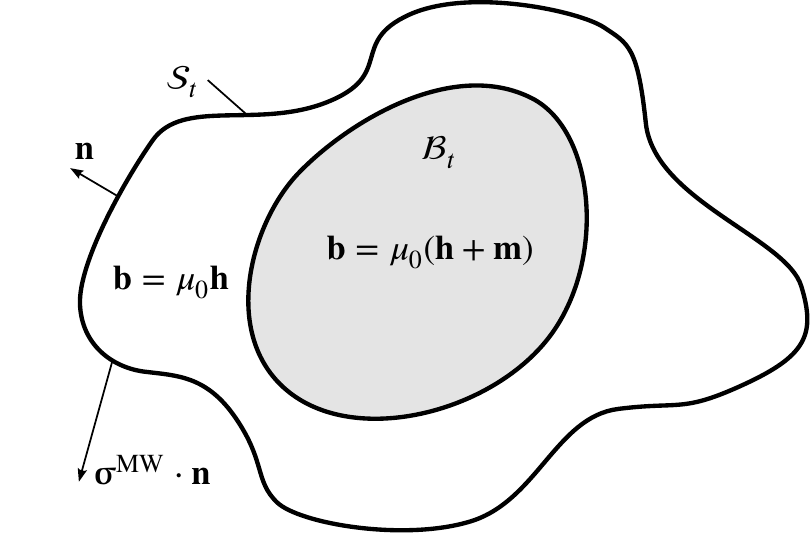}}%
    \caption{A magnetized body $\mathcal{B}_t$ surrounded by a non-magnetic medium. The closed surface $\mathcal{S}_t$ lies in the non-magnetic domain and encloses the magnetic body.}
    \label{fig:magnetized_body_in_free_space_surface}
\end{figure}

When considering two rigid magnetic bodies that are separated by a non-magnetic
medium, there will be some kind of force transmitted
between these two bodies. The transmission indeed happens through the
non-magnetic medium. In fact, the net magnetic force $\mathcal{F}^\text{mag}$ on a magnetized body can
be computed by the surface integral of the normal component of the Maxwell
stress over any closed surface $\mathcal{S}$ around the body of interest
\citep{brown1966}
\begin{align}
    \label{eq:net_force_on_body}
    \mathcal{F}^\text{mag} = \int_{\mathcal{S}_t} \Bsigmamw \cdot \bn \da .
\end{align}
See Figure \ref{fig:magnetized_body_in_free_space_surface} for an illustration.

Orthogonal to that is the important case of a magnetized body that does not experience
a net magnetic force. This can be a body with remanent magnetization that is remote from
other magnetic media or a uniformly magnetized body in a uniform external field.
Then, even in absence of a net force, there 
will still be jumps in the normal component of the magneto-mechanical stress across the
boundary of the magnetic body that depend on the magnetization $\bbm$ (see equations
\eqref{eq:magnetic_tractions_b} and \eqref{eq:magnetic_tractions_b}). 
In other words, the body still experiences tractions that are equilibrated in the sense
that the net force vanishes. However, they may still cause deformation.
The overall effect of these traction on a body depends on its shape, which is
the fundamental mechanism behind the shape-dependent magneto-mechanical
response of MREs \citep{keip+rambausek2017}. A corresponding energetic
point of view has been put forward by \citet{rambausek+keip2018}.
Such settings are also formidable test cases for the present study because of
the strong coupling of boundary deformation and magnetic response even in
geometrically ``simple'' settings.

%We close this subsection with the remark that this work is not concerned with
%detailed modeling of magnetic materials, which is a very active topic on its
%own. Instead we focus on the effect of their presence on the numerics of
%surrounding non-magnetic media, which as outlined above play a crucial role as
%transmitter of magnetic fields and forces.

\section{Spurious magneto-mechanical coupling}
\label{sec:spurious_forces}

\begin{figure}[!ht]
    \centering
    \ifthenelse{\equal{\figuremode}{\fmodesvg}}{%
    \includesvg[width=\textwidth,pretex=%
    \small%
    \newcommand{\labela}{\normalsize(a) boundary value problem}%
    \newcommand{\labelb}{\normalsize(b) unphysical defo. (energy form.)}%
    \newcommand{\labelc}{\normalsize(c) unphysical defo. (energy--co-energy form.)}%
    \newcommand{\lbext}{\normalsize$\bb^\infty$}%
    \newcommand{\llength}{\normalsize$L=20R$}%
    \newcommand{\lrad}{\normalsize$R$}%
    \newcommand{\lmagn}{\normalsize{$\volb^\text{magn}$ (magnetic)}}%
    \newcommand{\lnonmagn}{\normalsize{$\volb^\text{nonm}$ (non-magnetic)}}%
    \newcommand{\lbounfree}{\normalsize{$\bounfree$}}%
    ]{spurious_forces_example_circle_new}}{%
    \includegraphics{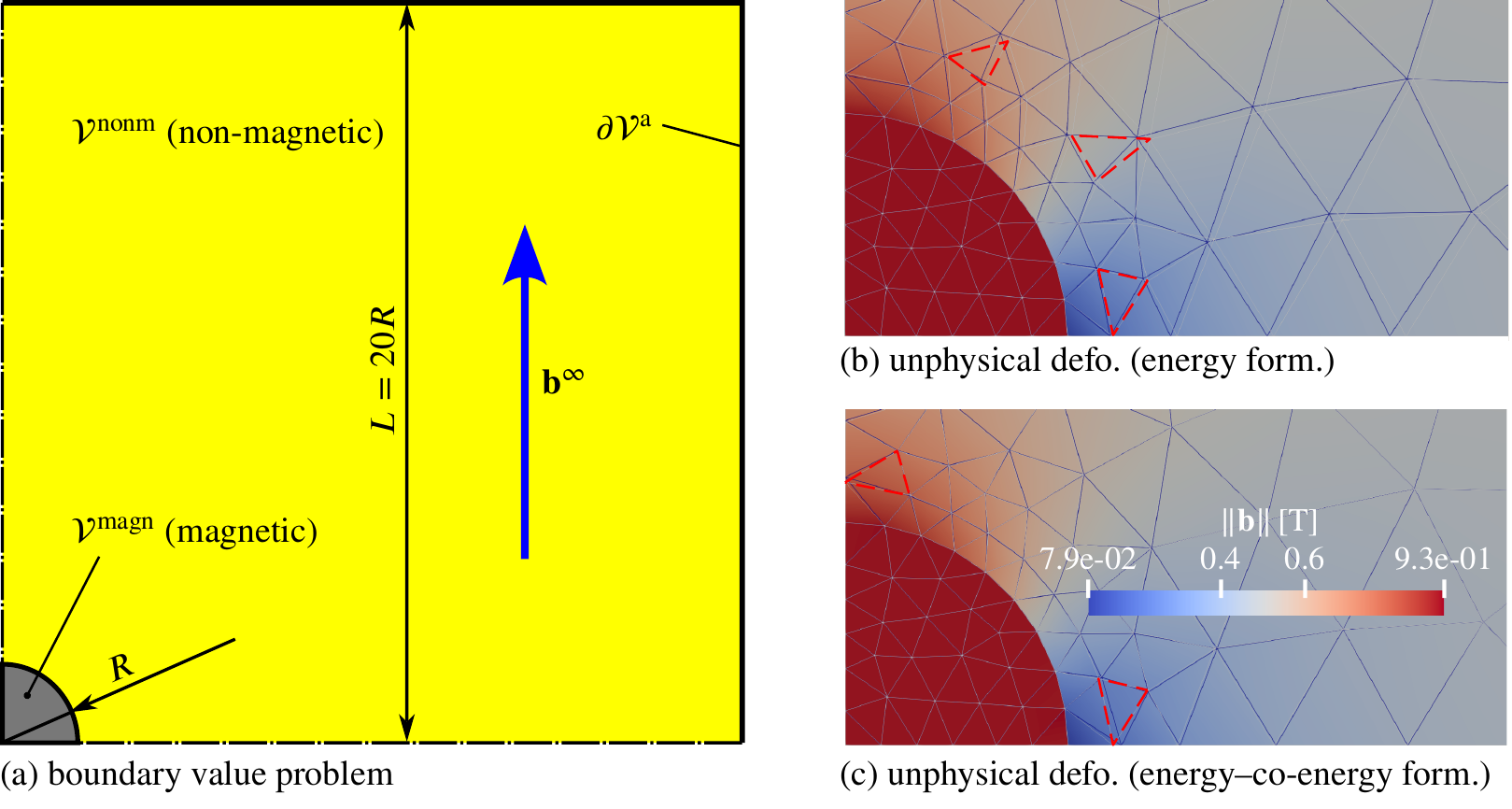}}%
    \caption{Test problem demonstrating a mild case of spurious deformation in a   
        non-magnetic deformable elastic domain.
        Subplot (a) depicts the (quarter) boundary value problem of a rigid
        magnetic particle $\volb^\text{magn}$ embedded in a non-magnetic elastic medium
        $\volb^\text{nonm}$ exposed to a
        uniform external magnetic field $\bb^\infty$.
        Subplot (b) show the result obtained with the magneto-mechanical energy
        formulation \eqref{eq:magnetomechanical_energy_principle} with color-contours 
        indicating to the magnetic field magnitude.
        The blue mesh corresponds to the deformed configurations, whereas
        the grey lines indicated the undeformed mesh.
        The dashed red triangles highlight parts of the \emph{undeformed} configuration
        which are significantly deformed under applied field.
        The corresponding plot for the energy--co-energy formulation
        \eqref{eq:magnetomechanical_energy_co-energy_principle} is depicted in subplot
        (c).
    }
    \label{fig:spurious_forces_example}
\end{figure}

While the theory clearly states that non-magnetic media do not deform in
reaction to magnetic fields, one may observe the
contrary in numerical simulations.
For the purpose of a brief demonstration, we consider the magneto-mechanical
(quarter) boundary value problem (BVP) of a circular quasi-rigid magnetic body 
surrounded by a rather soft non-magnetic medium such as a very soft elastomer or some
extremely soft auxiliary material representing air.
The geometrical setting is illustrated by
Figure~\ref{fig:spurious_forces_example}a.
The BVP is formulated in terms of both coupled variational principles
\eqref{eq:magnetomechanical_energy_principle} and
\eqref{eq:magnetomechanical_energy_co-energy_principle}, which are discretized
with second-order isoparametric Lagrange elements for both the deformation (displacements) and the out-of plane vector
potential component or the
scalar magnetic potential, respectively, as commonly done in finite magneto- and
electro-elasticity
\citep{%
    vu+steinmann+possart2007,%
    javili+chatzigeorgiou+steinmann2013,%
    keip+steinmann+schroeder2014,%
    miehe+vallicotti+teichtmeister2016,keip+rambausek2016,pelteret+etal2016,danas2017}.
The resulting fully coupled monolithic
nonlinear system obtained is solved by a Newton-Raphson procedure.
Before we continue, we need to clarify some terms. With ``fully coupled'' we 
simply refer to schemes which do not neglect any coupling that comes out of the
variational principle to be discretized, i.e. no physical effects present in
the theory are neglected or dropped for computational considerations or
convenience. Second, ``monolithic'' means that one solves the equations
representing the fully coupled discretized system for all primary variables at
the same time in the whole discretized spatial domain. 
For example, in the
case of \eqref{eq:magnetomechanical_energy_co-energy_principle}, one solves for
the pair $\{\defomap, \phi\}$ in all magnetic and non-magnetic media under
consideration.

Below we consider both the energy ($\{\defomap,\bA\}$) and energy--co-energy 
formulations ($\{\defomap,\phi\}$).
For both the magnetic and non-magnetic domain we consider the energy densities
\begin{align}
    \label{eq:demo_magnetic_material_B}
\PsiCB = \frac{G}{2} \left[\Tr\bC - 2\log J - 3\right] + 
\frac{G'}{2} \left(J - 1\right)^2 - 
\frac{\chi}{2 J \mu_0(1+\chi)} \bB\cdot(\bC\cdot\bB)
+ \PsiCBvac
\siand
    \label{eq:demo_magnetic_material}
    \PsiCH = \frac{G}{2} \left[\Tr\bC - 2\log J - 3\right] + 
    \frac{G'}{2} \left(J - 1\right)^2 - 
    \frac{J \mu_0 \chi}{2} \bH\cdot(\bC^{-1}\cdot\bH)
    + \PsiCHvac,
\end{align}
where vacuum permeability is
$\mu_0=\SI{0.4\pi}{\micro\tesla\meter\per\ampere}$.
For either formulation, the material parameters employed for the practically
rigid magnetic domain are 
$\{G,G',\chi\} = \{\SI{1e4}{\kilo\pascal},\SI{5e5}{\kilo\pascal},10\}$. 
For the non-magnetic medium we choose 
$\{G,G',\chi\} = \{\SI{1e-5}{\kilo\pascal},\SI{5e-4}{\kilo\pascal},0\}$ 
in when using the energy formulation but
$\{G,G',\chi\} = \{\SI{1e-3}{\kilo\pascal},\SI{5e-2}{\kilo\pascal},0\}$ in
the energy--co-energy formulation.
For the non-magnetic domain, we for the purpose of demonstration employ
different mechanical parameters because the spurious coupling
is more pronounced for the energy--co-energy formulation in this example.

%The magnetic loading is applied via the ``external'' potential
%\begin{align}
%    \label{eq:magnetic_loading_potential}
%    \Piext = \int_{\bounf_0} -(\bb^\infty \cdot \bN)\,\phi\dA,
%\end{align}
%where $\bN$ is the outward-pointing unit normal to the surface. The corresponding
%loading term in the weak form reads
%\begin{align}
%\label{eq:magnetic_loading_weak_form}
%\delta\Piext = \int_{\bounf_0} -(\bb^\infty \cdot \bN)\,\delta\phi\dA .
%\end{align}

\begin{figure}[!ht]
    \centering
    \ifthenelse{\equal{\figuremode}{\fmodesvg}}{%
    \includesvg[width=\textwidth,pretex=%
    \newcommand{\labela}{(a)}% magnetostatic co-energy}%
    \newcommand{\labelb}{(b)}% deviation from rigid case}%
    \newcommand{\lbext}{$b^\infty \perunit{\unitB}$}%
    \newcommand{\lhform}{energy--co-energy ($h$) form}%
    \newcommand{\lbform}{energy ($b$) form}%
    \newcommand{\lmagpot}{magnetostatic pot. $\Pi^{b|h} \perarbitraryunit$}%
    \newcommand{\lmagpotdiff}{$(\Pi^{b|h} - \Pi^{b|h}_{\text{rigid}}) \times 10^{6} \perarbitraryunit$}%
    ]{circle_in_air_energies_bh_with_grid}}{%
    \includegraphics{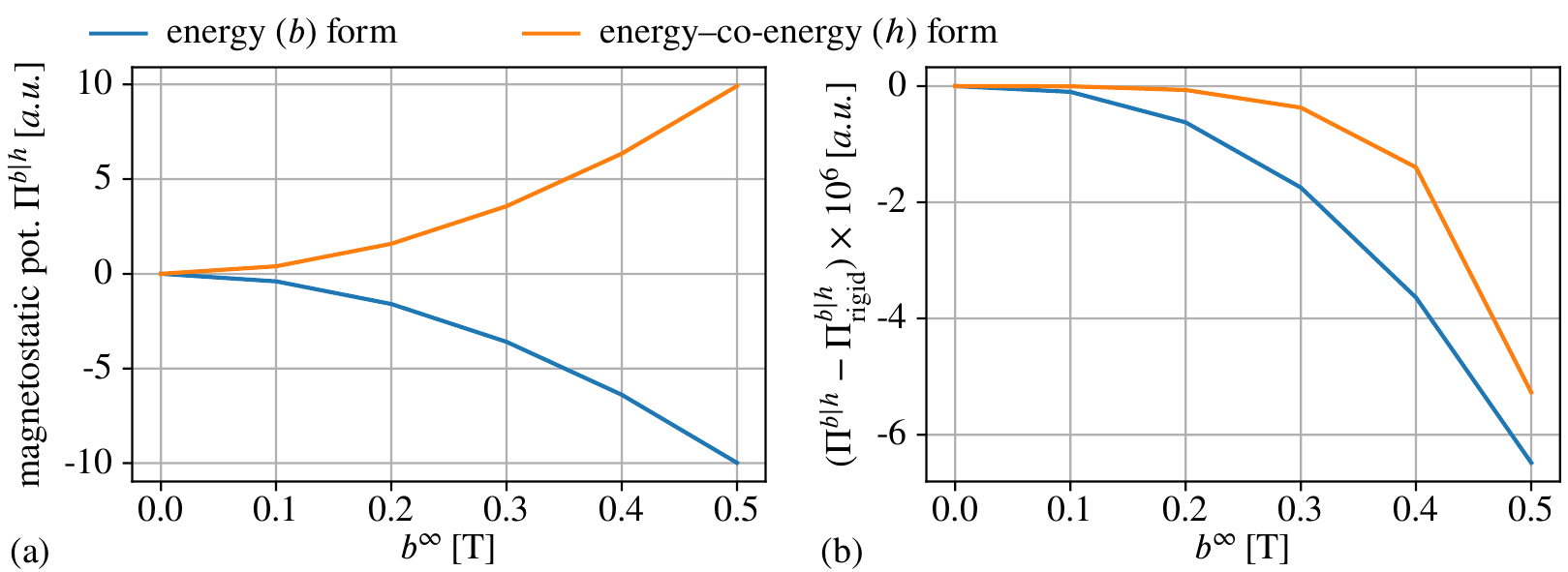}}%
    \caption{Magnetostatic (co-)energy $\Pi^{b|h}$ for a rigid magnetic particle
        embedded in a non-magnetic solid depending on the applied
        field magnitude $b^\infty$ (a). 
        The respective difference to the (co-)energy $\Pi^{b|h}_{\text{rigid}}$ obtained
        for a \emph{rigid} non-magnetic domain (b).}
    \label{fig:optimization_perspective}
\end{figure}

%\footnote{The results presented in this section can be reproduced with the
%    Jupyter Notebook \texttt{spurious\_deformation\_naive\_circle\_demo.ipynb},
%which is provided as supplementary material.
%For readability, detailed input data, e.g., material models and parameters, are
%omitted for the sake of readability of this explanatory section.}
Figure~\ref{fig:spurious_forces_example}(b) and (c) show color-contours 
of the magnetic field magnitude in the deformed configuration obtained with
the magneto-mechanical energy formulation \eqref{eq:magnetomechanical_energy_principle}
and the energy--co-energy formulation
\eqref{eq:magnetomechanical_energy_co-energy_principle}, respectively.
The blue mesh corresponds to the deformed configuration while the light grey mesh
is the undeformed one. Both meshes coincide in the practically rigid magnetic particle.
From the physical theory one would not expect any
deformation or displacements in the non-magnetic domain, since the magnetic body is
spatially fixed and rigid whereas the deformable material is non-magnetic and clamped at the outer boundary $\bounfree$. 
However, the dashed-red triangles marking undeformed triangle cells which undergo
significant deformation (dashed red versus solid blue) show that the opposite is the
case. This clearly is against the physical theory and, as we argue below, a universal
trait of fully coupled monolithic discretizations.
The fundamental problem with such approaches is that the solution fields
differ for each discretization, which is an inherent property of any numerical
method in general. This means that, 
the values of the governing energies or potentials evaluated for numerical
solution states change with the discretization. 
The crucial step is now to
consider the magnetostatic (co-)potentials
\begin{align}
    \Pi^B(\defomap, \bA) &= \int_{\volf_0} \PsiCBvac(\bF, \bB) \dV
     - \int_{\volb^\text{magn}_0} \frac{\chi}{2 J \mu_0(1+\chi)} \bB\cdot(\bC\cdot\bB) \dV
     - \int_{\bounf_0} (\bh^\infty\times\bn)\cdot\bA \dA
     \nonumber \\
     &= \int_{\volf_t} \psibvac(\bb) \dv 
     + \int_{\volb^\text{magn}_t} \frac{\chi}{2 \mu_0(1+\chi)} \norm{\bb}^2 \dv
     - \int_{\bounf_t} (\bh^\infty\times\bn)\cdot\ba \da = \Pi^b(\ba;\defomap)
%    & \text{and} &
\siand
    \Pi^H(\defomap, \phi) &= \int_{\volf_0} \PsiCHvac(\bF, \bH) \dV 
    - \int_{\volb^\text{magn}_0} \frac{J \mu_0 \chi}{2} \bH\cdot(\bC^{-1}\cdot\bH)
    - \int_{\bounf_0} (\bB^\infty\cdot\bN)\phi \dA
    \nonumber \\
    &= \int_{\volf_t} \psih(\bh) \dv 
       - \int_{\volb^\text{magn}_t} \frac{\mu_0 \chi}{2} \norm{\bh}^2 \dv
       - \int_{\bounf_t} (\bb^\infty\cdot\bn)\phi \da = \Pi^h(\phi;\defomap)
    ,
\end{align}
respectively, where $\defomap$ is considered as a ``given'' parameter in $\Pi^b$ and
$\Pi^h$.
Recall from Sections~\ref{sec:intro_magnetostatics} and
Sections~\ref{sec:intro_magnetoelasticity} that -- for a
given $\defomap$ -- the magnetostatic problem is solved
for $\ba$ by minimizing $\Pi^b$ and
for $\phi$ by maximizing $\Pi^h$, respectively. However, we also solve for $\defomap$
via the minimization of the total (co-)energy%
\footnote{We use the short $\Pi^{B|H}$ when
    we refer to both $\Pi^{B}$ and/or $\Pi^{B}$.} 
$\Pi^{B|H}$. When we were not in possession of numerical but exact magnetostatic
solutions, $\Pi^{B|H}$ would be invariant with respect to $\defomap$ in the nonmagnetic
domain (the magnetic domain is considered to be rigid in the present example). 
But as our numerical solutions are not exact, some deformation has to happen%
\footnote{The deformation magnitude, of course, depends on the stiffness of the
    non-magnetic medium.} as we have seen in Figure~\ref{fig:spurious_forces_example}b,c.
The energetic considerations are illustrated by
Figure~\ref{fig:optimization_perspective}, 
where subplot (a) shows the magnetostatic potentials $\Pi^{b|h}$ at a solution state
and (b) depicts the difference to the rigid case.
The graphs in Figure~\ref{fig:optimization_perspective}a are practically mirrored
(neglecting FE errors) 
as is expected from the duality of the energy and the co-energy formulations in the case
of linear magnetic materials.
In Figure~\ref{fig:optimization_perspective}b we observe that both potentials
$\Pi^{b}$ and $\Pi^h$ are smaller for the deformable non-magnetic medium than for the
rigid case, i.e. there was a minimization of the magnetostatic energy through
deformation.
In the energy formulation this is \textit{formally} favorable as
we interprete a smaller computed total energy as being closer 
to the (true) solution than a higher total energy. But of course, 
this minimization is only possible because of numerical errors and the resulting
deformation is entirely unphysical and undesirable.
In the energy--co-energy case, we do not have any ``formally'' favorable outcome.
In particular, the deformation pattern is to our experience in general quite different
from the one obtained for the energy formulation. Moreover, it is typically more
pronounced -- recall that we had $G=\SI{1e-3}{\kilo\pascal}$ in the energy--co-energy
formulation but $G=\SI{1e-5}{\kilo\pascal}$ in the energy formulation.
In the present example, the energy--co-energy would crash for applied field magnitudes
somewhere between around $b^\infty=\SI{0.5}{\tesla}$ and $b^\infty=\SI{0.6}{\tesla}$.
In contrast, the energy formulation runs through until at least
$b^\infty=\SI{1.0}{\tesla}$.
While the latter appears to be more robust, it might also crash eventually.
However, if the soft non-magnetic medium is indeed a physical material and not
some auxiliary model for an air-like environment, any spurious deformation must
be avoided. 
In the next section we present two approaches that achieve this goal.

\begin{rmk}
    \label{rmk:reference_to_r-adaptivity}
    The above considerations have parallels to those behind \emph{variational
    r-adaptivity} for finite elasticity
    \citep{thoutireddy+ortiz2004,kuhl+askes+steinmann2004,askes+kuhl+steinmann2004}.
    In that context, the ``spatial motion'' (the ``usual'') finite strain
    elasticity problem on the Lagrangian configuration plays a role comparable
    to the magnetostatic problem in the Eulerian configuration.
    The Lagrangian mesh optimization, on the other hand, is governed
    by the ``material motion'' problem which is parallel to the ``spatial''
    motion problem -- in a sense \textit{Eulerian mesh optimization} -- in the
    present context.
    We furthermore point out that the ``material motion'' problem corresponds
    to configurational mechanics and the respective forces are
    \emph{configurational forces} \citep{gurtin1994,gurtin1999,eshelby1999}. 
    They are expected to
    vanish everywhere except at material \emph{defects} such that one may attain
    the viewpoint that the discretization of a body indeed introduces defects.
\end{rmk}

\section{Two novel approaches to cure spurious magneto-mechanical coupling in
         non-magnetic media}
\label{sec:remedies}

Without loss of generality, we from now on restrict ourselves 
to the use of the magnetostatic co-energy and the variational principle
\eqref{eq:magnetomechanical_energy_co-energy_principle}. 
Nevertheless, we emphasize that
all derivations and observations presented below have their equivalent counterpart
when \eqref{eq:magnetomechanical_energy_principle} is employed instead.
In particular, we recall that the Maxwell stress has a unique notion in non-magnetic
media. As a result, the methods introduced in this Section are directly applicable to
both the energy and energy--co-energy formulation.

Below we shall distinguish two cases. First, we consider magnetic bodies
surrounded by vacuum or air, i.e. a non-magnetic medium of practically negligible
stiffness. Then we do not have to
care too much about the actual deformation outside the bodies. 
What matters is that the deformation acceptable from a numerical perspective.
%By assumption, the deformation
%outside the bodies should not have any significant effect on the system anyways.
The second case is concerned with magnetic bodies that are embedded in a very
soft carrier medium such as an extremely soft elastomer, a gel or some
biological tissue. In that case one expects a specific mechanical
response of non-magnetic material such that we are not free to choose the mechanical
parameters. Because of this constraint, we regard the latter case as more difficult.

To our best knowledge, up to now, all remedies to spurious magneto-mechanical
coupling in non-magnetic media are concerned with spurious coupling in vacuum
or air-like domains. Therefore, we start with this important case, where we
first briefly review existing approaches and then propose two novel cures to this
issue. They will be assessed by a comparison with each other and the staggered
scheme of \citet{pelteret+etal2016}. After that demonstration, in
Subsection~\ref{sec:remedies_for_solids} we will turn to
the extension of the proposed methods to extremely soft, non-magnetic solids.

Before jumping into the details, we remark that we only consider the case when
spurious deformations are an issue only in the non-magnetic domains but not in the
magnetic ones. In other words, throughout this work we assume magnetic domains to be
sufficiently stiff such that spurious deformations are negligible.
We anticipate though, that this might not always be the case.

\subsection{Eliminating spurious coupling in vacuum and air-like domains}
\label{sec:remedies_for_vacuum}

In case of vacuum or air surrounding the body of interest, one might look into schemes
which just do not solve the coupled problem in a monolithic but in a staggered way
along the lines of \citet{pelteret+etal2016}. In their scheme, one first solves
the fully coupled problem in a monolithic may but with deformation blocked in the
interior of the non-magnetic domain. By that one ends up with a thin deformable
layer between the boundary of the solid bodies and the first FE nodes in the
vacuum domain with vanishing stiffness.
As a consequence, this first part of problem is free of any artificial elastic
model for the air or vacuum domain.
The second step consists of smoothing the displacement field in the vacuum
domain by an auxiliary problem. 
While this adaption does not have any effect on the stiffness experienced by
the solid bodies from their ``empty'' surrounding, there is an effect on the
magnetostatic solution due to the change of the Eulerian configuration. And
since the problem is physically coupled, the mesh adaption step perturbs the
equilibrium computed in the first, coupled step. This can be a disadvantage
when very soft bodies or compliant (slender) structures are considered.
The situation can be improved by iterating between the coupled and the 
displacement-smoothing step.

A second highly successful scheme for the treatment of vacuum or air in this
context is the use of non-local algebraic constraints that bind the displacement
in the vacuum domain to the motion of the boundaries of the embedded bodies
\citep{psarra+etal2019}. This scheme is in fact monolithic and thus does not
suffer from the convergence issues of the staggered approach discussed before.
However, the non-local constraints can be difficult to setup for complex
geometries \citep{dorn+bodelot+danas2021} and furthermore lead to an
unfavorable sparsity structure of the system matrix.
Both, the staggered and the non-local constraint approaches -- if applicable --
lead to quantitatively better results than the probably most straightforward
method from an implementation perspective:
assigning a small but sufficient artificial stiffness to the vacuum domain
\citep{keip+rambausek2016,dorn+bodelot+danas2021}. By that one ends up with a
perturbed boundary value problem similar to the one employed in
Section~\ref{sec:spurious_forces} for the demonstration of spurious magnetic
forces. We will refer to this approach as ``naive monolithic scheme''.
As we have seen in Section~\ref{sec:spurious_forces},
spurious coupling may eventually lead to crashes.
Thus, for the ``naive monolithic scheme'' we face the problem that the spurious
coupling sets a lower bound on the (auxiliary) stiffness added.
Since the accuracy of the ``vacuum'' approximation depends on the ratio of
``solid'' to ``vacuum'' stiffness, the lower bound on the auxiliary stiffness
will at some point lead to significant deviations from the unperturbed problem
and its solution. The ``naive'' approach can be significantly improved by adapting the
stiffness added to the magnitude and the ``profile'' of the magnetic fields in
the vacuum domain as demonstrated by \citet{rambausek+mukherjee+danas2022}. 
However, this requires a certain knowledge of the problem or
complicated heuristics to automatically assign optimal stiffness parameters.
Moreover, such improvements only shift the limit but do not solve the
underlying problem.

A somewhat distinct method is that of
\citet{liu+etal2020} which is essentially based on magnetic forces and
tractions derived from
numerically computed fields and (localized and free) electric currents. By
that, no forces emerge in non-magnetic domains. In order to adapt the deformation 
of the vacuum domain to that of the
slender structure under consideration they thus may assign a much smaller --
indeed negligible -- artificial stiffness.
However, their specific scheme for the computation of magnetic
volume and surface force densities requires to compute derivatives of the
magnetization, which are not directly available in standard FE simulations 
and thus prevents a monolithic approach.
Instead, their solution procedure is a staggered one based on the alternate
solution of the magnetostatic and the elasticity problem. Hence, one needs a certain
number of iterations between the two subproblems in order to solve the coupled
problem.

What we propose below has partially great similarities with the method of
\citet{liu+etal2020} on a higher level in the sense that we focus on what is going on
at the interface between magnetic and non-magnetic bodies. However, in contrast to
\citep{liu+etal2020}, our schemes essentially rely on the same building blocks as the
naive monolithic scheme. By that the proposed methods can be implemented with minimal
adaptions to existing code. In particular, there is no need
to compute forces emerging from currents nor anything else that would require a
staggered algorithm. Thus, the new schemes can be easily
linearized and will be implemented in a monolithic way. When starting from an
implementation of the naive monolithic scheme, only minor adaptions of code are
required.

\subsection{The Maxwell traction approach}
\label{sec:maxwell_traction_approach}

For the derivation of the first approach we consider the ``mechanical'' weak form
corresponding to the first variation of the governing potentials with respect to
deformation restricted to vacuum domains $\volfree$.
We emphasize that at this stage, the only (co-)energy density considered in vacuum
is the magnetostatic one. Thus, the only stress present is the vacuum Maxwell
stress \eqref{eq:cauchy-type_maxwell_stress_in_vacuum} 
and the ``mechanical'' weak forms for the energy and the co-energy versions
coincide
\begin{align}
    \label{eq:magnetostatic_weak_form_in_vacuum_1}
    \delta_{\defomap} \int_{\volfree_0} \PsiCBvac(\bC, \bB) \dV =
    \delta_{\defomap} \int_{\volfree_0} \PsiCHvac(\bC, \bH) \dV =
    \int_{\volfree_t} \Bsigmamw \mcolon \grad{\delta{\defomap}} \dv .
\end{align}
As discussed in the context of \eqref{eq:cauchy-type_maxwell_stress_in_vacuum},
$\Bsigmamw$, the vacuum Maxwell stress is a divergence-free field such that we
may safely%
\footnote{``Safe'' in the sense that the concerns raised by
     \citet{reich+rickert+mueller2017}, who warn from confusing long and short range
     forces, do not apply.}
transform the equation above by means of the divergence theorem
\begin{align}
\label{eq:magnetostatic_weak_form_in_vacuum_2}
 \int_{\volfree_t} \Bsigmamw \mcolon \grad{\delta{\defomap}} \dv =
 \int_{\bounfree_t} \delta{\defomap} \cdot \left(\Bsigmamw \cdot \bn\right) \da =
 \int_{\bounfree_0} \delta{\defomap} \cdot \left(\PKonemw \cdot \bN\right) \dA,
\end{align}
where $\PKonemw$ is the first-Piola-Kirchhoff-type instance of the Maxwell
stress. The important observation here is that, in the case the FE discretizations
under consideration\footnote{In this work we are concerned with Lagrange-type or
    equivalent elements for the deformation map or displacements, respectively.}, 
the boundary integrals in \eqref{eq:magnetostatic_weak_form_in_vacuum_2} do contribute
any discrete ``interior'' forces\footnote{``Interior'' in the sense that the
    correspond to degrees of freedom that control the deformation in the
    \emph{interior} of the nonmagnetic domain.}.
As a consequence, no spurious coupling within $\volfree$ can emerge from such boundary
integrals whereas the volume integral in
\eqref{eq:magnetostatic_weak_form_in_vacuum_2} indeed suffers from this numerical
artifact.%
\footnote{There is still some spurious coupling within the interface or
    boundary $\bounfree$, but it is assumed that the medium outside is stiff enough
    such that spurious deformation is suppressed. If this is not the case, the only
    remedy is sufficiently accurate magnetostatic solution.}
This is somehow parallel to the use of magnetic surface forces by \citet{liu+etal2020}
but approached from a quite different angle. Indeed, the jump of the normal components
of the ``magnetic'' stresses corresponds to the magnetic surface tractions (see
\eqref{eq:magnetic_stresses_b} and \eqref{eq:magnetic_stresses_b}), but here we
only compute the respective contribution on the vacuum side.

In order to obtain a non-singular system matrix in FE simulations
we equip the vacuum domain with a very soft elastic material. However, 
in contrast to the ``naive'' monolithic scheme, the artificial stiffness does
not have to counterbalance the effects of spurious forces and is thus only
limited from below by requirements of the linear solver and numerical stability,
respectively.

Summarizingly, the above considerations ask us to implement the right-hand-side of
\begin{align}
    \label{eq:maxwell_traction_problem}
    \delta_{\defomap} \int_{\volfree_0} \PsiCHvac(\bC, \bH) + \Psi^\text{a}(\bC) \dV =
    \delta_{\defomap} \int_{\volfree_0} \Psi^\text{a}(\bC) \dV
    + \int_{\bounfree_0} \delta{\defomap} \cdot \left(\PKonemw \cdot \bN\right) \dA
\end{align}
instead of the left-hand-side,
where $\Psi{^\text{a}}$ is the auxiliary strain energy density.
Finite element implementations of the boundary integral in
\eqref{eq:maxwell_traction_problem}
approach are briefly outline below. The actual implementations in
\texttt{Netgen/NGSolve} are provided as supplementary material.

\paragraph{Direct implementation of the boundary integral}
%\label{sec:mw_traction_impl_boundary_integral}
The most obvious finite element implementation of the ``Maxwell traction''
approach is the implementation of the boundary integral in
\eqref{eq:maxwell_traction_problem} just as what it is -- a boundary integral. 
The essential additional algorithmic ingredients are 
\begin{itemize}
    \item a routine the gathers all element faces that belong to the boundary
    of the non-magnetic domain and of which the parent (volume) element belongs
    to the vacuum domain and
    \item the evaluation of the outward-pointing unit normal on these faces.
\end{itemize}
For our implementation we choose \texttt{Netgen/NGSolve} which directly
provides the outward-point unit normal vector and a symbolic language for the
integration over element boundaries.
However, as will be shown in the numerical examples in
Section~\ref{sec:comparison_in_air}, the convergence for
refined meshes of such an implementation is not satisfactory.
We also point out that while the boundary integral can readily be linearized to obtain
its contribution to the system matrix, this contribution is non-symmetric.

\paragraph{Implementation via discrete force omission}
Alternative to the boundary integral implementation, one can exploit the fact
that $\PKonemw$ is expected to be divergence free in a non-magnetic domain,
such that all resulting discrete ``interior'' forces can just be ignored.
In contrast, the discrete forces corresponding to boundary deformation degrees of
freedom balance discrete forces stemming from applied tractions
or the body on the other side of an interface.
Hence, instead of actually computing the boundary integral in strict sense, we may
simply compute the volume integral in \eqref{eq:maxwell_traction_problem}.
This is what one
would also do in the ``naive'' monolithic approach but ignoring all discrete forces
resulting from $\delta\PsiCHvac=\PKonemw:\delta\defomap$ that correspond to
deformation degrees of freedom (DoFs)
in the \emph{interior} of $\volfree$.
The actual implementation thus only needs to be able
to distinguish between discrete forces in the interior of the vacuum domain and its boundary. Depending on the finite element software, the deletion of
forces might be done by simply omitting the respective contribution
during assembly. Alternatively, one can assemble the residual vectors resulting from
$\Psi^\text{a}$ and $\PsiCHvac$ separately in a first step. In a second step,
the ``interior'' forces from the latter residual vector are deleted such that only
the boundary forces remain. The final residual vector is the obtained as the sum
of the auxiliary residual vector and the fleshed-out vacuum residual.
It is important to note that the omission or deletion of discrete force
contribution entails the deletion of the corresponding contributions to the
system matrix, rendering the latter non-symmetric.

Our actual implementation employs the latter method as this can be done without
any low-level changes to the assembly routines of \texttt{NGSolve}.

\subsection{The traction compensation approach for air-like media}
\label{sec:traction_compensation_approach}

The second approach proposed also relies on tractions but, besides that, builds
upon a converse line of thinking. 
The basic idea is to assign an auxiliary finite strain
energy $\Psi^\text{a}$ to the 
vacuum domain then compensate the effect of the stiffness added by subtracting
the corresponding \emph{mechanical} traction from the boundary $\bounfree$.
This means adding an ``integral zero'' to the functional
under the \emph{condition} that the corresponding stress field is \emph{divergence-
free}.
We start from the weak form \eqref{eq:magnetostatic_weak_form_in_vacuum_1}
and proceed as outlined to obtain
\begin{align}
    \label{eq:stiffness_compensation_1}
    &\int_{\volfree_t} \Bsigmamw \mcolon \grad{\delta{\defomap}} \dv
    = 
    \int_{\volfree_t} \Bsigmamw \mcolon \grad{\delta{\defomap}} \dv
    + \underbrace{
        \int_{\volfree_t} \Bsigma^{\text{a}} \mcolon \grad{\delta{\defomap}} \dv - 
        \int_{\bounfree_t} \delta{\defomap} \cdot \left(\Bsigma^{\text{a}} \cdot \bn \right) \da
    }_{ = - \int_{\volfree_t} \delta{\defomap} \cdot \left(\div \Bsigma^{\text{a}}\right) \dv = 0}
\sitext{with}
\label{eq:def_sigma_aux}
   &\Bsigma^{\text{a}} = J^{-1} \pdiff{\Psi^\text{a}}{\bF} \cdot \bFT .
\end{align}
While the right-hand-side of \eqref{eq:stiffness_compensation_1} is the basis for our
implementations,
further manipulations provide additional insight into the method:
\begin{align}
\label{eq:stiffness_compensation_2}
&\int_{\volfree_t} \Bsigmamw \mcolon \grad{\delta{\defomap}} \dv +
\int_{\volfree_t} \Bsigma^{\text{a}} \mcolon \grad{\delta{\defomap}} \dv - 
\int_{\bounfree_t} \delta{\defomap} \cdot \left(\Bsigma^{\text{a}} \cdot \bn \right) \da
\nonumber\\
&=
\int_{\volfree_t} 
\underbrace{\Bsigma}_{\Bsigmamw + \Bsigma^{\text{a}}} 
    \mcolon \grad{\delta{\defomap}} \dv - 
\int_{\bounfree_t} 
\delta{\defomap} \cdot \left(\Bsigma^{\text{a}} \cdot \bn \right) \da
\nonumber\\
&=
\int_{\volfree_t} 
-\delta\defomap \cdot \left(\div\Bsigma\right) \dv
+ 
\int_{\bounfree_t} 
\delta{\defomap} \cdot \left(\Bsigmamw \cdot \bn \right) \da .
\end{align}
From this equation one can see that, irrespective of the actual value of
$\Bsigma^{\text{a}}$, only the Maxwell stress is exerted on embedded bodies
(``across'' the boundaries of the vacuum domain). Thus, as the actual
deformation in the \emph{interior} of the vacuum domain does not affect its
neighboring bodies, which allows for great freedom in the choice (form and
magnitude) of $\Psi^\text{a}$ or $\Bsigma^\text{a}$, respectively.
A second aspect is that one effectively solves for 
$\div\Bsigma=\div(\Bsigmamw+\Bsigma^\text{a})=0$ which allows for 
$\Bsigma^\text{a} \neq 0$ such that the auxiliary stiffness can counterbalance
any spurious forces corresponding to $\Bsigmamw \neq 0$.
In that case, the terms added in the derivation do not exactly add up to zero.

The crucial aspect from an implementation point of view is the treatment of the
boundary integral in \eqref{eq:stiffness_compensation_1}.
In the remainder of this subsection we present two possible finite element
approaches for which the actual implementation in \texttt{Netgen/NGSolve} is
provided as supplementary material.

\paragraph{Direct computation of the traction compensation integral}
Similar to the Maxwell traction approach, we start with the direct
implementation of \eqref{eq:stiffness_compensation_1}. 
The volume integrals are
exactly the same as for the ``naive'' monolithic implementation. The only
addition in this sense is the boundary integral in that equation, of which the
implementation is technically parallel to that for the Maxwell tractions.
It also has the same drawback of bad convergence under mesh refinement as will
be shown in the numerical examples in Section~\ref{sec:comparison_in_air}. Similarly, the resulting system matrix is non-symmetric.

\paragraph{Implementation via discrete force omission}
Reusing the viewpoint of force omission approach for Maxwell tractions,
one can implement traction compensation via manipulations of the discrete force
vector resulting from standard evaluation of volume integrals and the corresponding
system matrix contributions.
However, in contrast to the Maxwell tractions, now one leaves the discrete
forces corresponding to displacements in the interior as they are but removes
the auxiliary discrete forces that correspond to
displacement DoFs at the boundary of the vacuum domain. 
This means the we may think of implementing \eqref{eq:stiffness_compensation_1} as 
\begin{align}
    \label{eq:stiffness_compensation_1b}
&\int_{\volfree_t} \Bsigmamw \mcolon \grad{\delta{\defomap}} \dv
= 
\int_{\volfree_t} \Bsigmamw \mcolon \grad{\delta{\defomap}} \dv
+   \int_{\volfree_t} \Bsigma^{\text{a}} \mcolon \grad{\delta{\defomap}} \dv - 
\underbrace{
    \int_{\volfree_t} \Bsigma^{\text{a}} \mcolon \grad{\delta{\defomap}} \dv
}_{\text{omit ``interior'' forces}}
\shortintertext{or}
\label{eq:stiffness_compensation_1c}
&\int_{\volfree_t} \Bsigmamw \mcolon \grad{\delta{\defomap}} \dv
= 
\int_{\volfree_t} \Bsigmamw \mcolon \grad{\delta{\defomap}} \dv
 +
\underbrace{
    \int_{\volfree_t} \Bsigma^{\text{a}} \mcolon \grad{\delta{\defomap}} \dv
}_{\text{omit ``boundary'' forces}} .
\end{align}
Thus, for any implementation
the essential algorithmic ingredient is a mean to distinguish ``interior'' and
boundary'' displacement DoFs.
Our own implemention first assembles the volume integrals in 
\eqref{eq:stiffness_compensation_1c} separately, removes the boundary forces from 
the auxiliary residual vector and finally adds both vectors.
Accounting for the manipulations of discrete forces in the system matrix leads
to a loss of symmetry of the latter, as is the case for all approaches proposed
in this contribution.

\subsection{Assessment of the ``Maxwell traction'' and the ``traction compensation approach for air-like non-magnetic domains}
\label{sec:comparison_in_air}

We critically assess the proposed schemes by their performance in the boundary
value problem presented in Section~\ref{sec:spurious_forces}, in particular
Figure~\ref{fig:spurious_forces_example}a. The energy densities employed are of
the form \eqref{eq:demo_magnetic_material}. The material parameters for the
magnetic domain are $\{G, G', \chi\} = \{\SI{1}{\kPa}, \SI{50}\kPa, 10\}$.
We consider the non-magnetic domain to be air-like, where $\chi=0$ and $G=G'=0$.
However, depending on the numerical scheme applied, we employ different values for
the auxiliary hyperelastic parameters $G^\text{a}$ and ${G'}^\text{a}$. While the effect of their absolute value is to be studied, we keep the
ratio of $G^\text{a} / {G'}^\text{a} = 1/50$ throughout all examples.
Concerning the discretization we employ triangular meshes and 
second-order finite element spaces for each component of the deformation and for
the scalar magnetic potential. This appears to be a common choice in literature
such that the results reported in this work are directly relevant to research
applications.

From a theoretical perspective, the ``traction compensation'' approach is
exact whereas the ``Maxwell traction'' approach presented is not. 
This is because in the latter case the effect of the
auxiliary stiffness is very small but nevertheless finite.
In contrast, all effects from the auxiliary stiffness and body force
contributions are compensated in the former case.
However, despite being exact in the above sense, the results obtained with the
traction compensation approach still depend on the magnitude of added stiffness
and forces. This is particular observable for rather small added stiffness such
that the spurious coupling may still affect the deformation of the non-magnetic
domains. Therefore, in what follows, we investigate the effect of the absolute
value of the auxiliary shear modules $G^\text{a}$ and the discrete resolution
of the domain, i.e. the mesh size.

We start the comparison with individual parameter studies of the ``naive''
monolithic approach, continue with the ``Maxwell
traction'' and ``traction compensation'' implementations and finally present a
comparative convergence study.
All results presented below refer to the vertical displacement $u_2$ of the
point with initial position $A = (0, R)$, i.e. the intersection of the boundary of the
magnetic domain with the $y$-axis, at an applied magnetic field magnitude
$b^\infty=\SI{0.7}{\tesla}$. 

\begin{figure}[!ht]
    \centering
    \ifthenelse{\equal{\figuremode}{\fmodesvg}}{%
    \includesvg[width=\textwidth]{circle_meshes}}{%
    \includegraphics{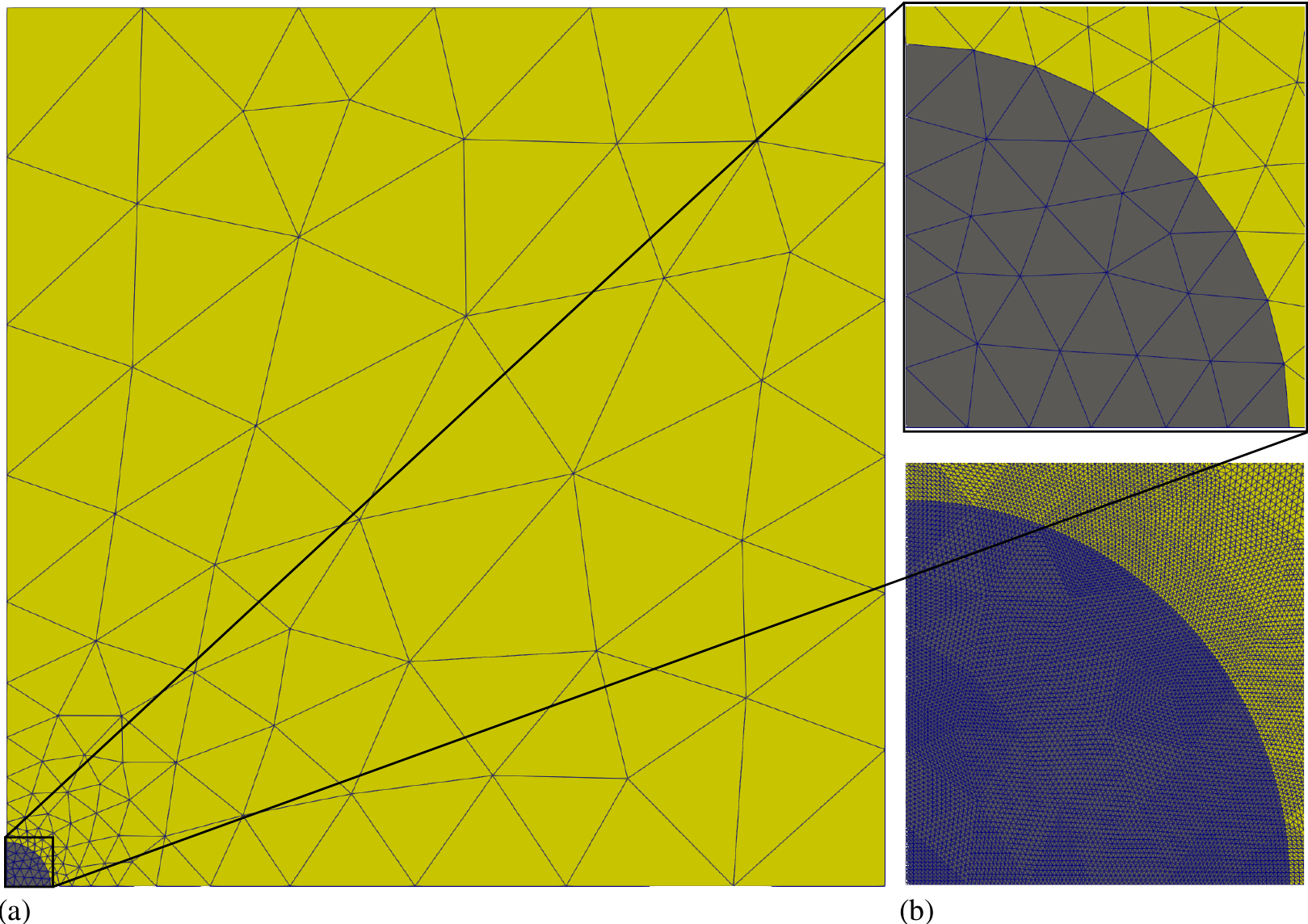}}%
    \caption{Coarsest and finest mesh of the circular magnetic domain embedded in a
         non-magnetic domain. Subplot (a) depicts the ``full'' (actually quarter
         domain) mesh. Details of the coarsest and the finest mesh employed are shown
         in (b). Note that actual geometry representation is of second order,
         which has only been lost during postprocessing.}
    \label{fig:circle_meshes}
\end{figure}

The coarsest mesh generated with \texttt{Netgen} is depicted in
Figure~\ref{fig:circle_meshes}a, a comparison of the
coarsest (refinement level zero) and the finest (refinement level four) mesh is
provided in Figure~\ref{fig:circle_meshes}b. Level zero corresponds to a total of
\num{222} cells and \num{1473} DoFs. After four refinement steps one has
\num{56832} cells and \num{343203} DoFs.

\subsubsection{The ``naive'' fully-coupled monolithic auxiliary stiffness approach}
\label{sec:circle_in_air_nv}

\begin{figure}[!ht]
    \centering
    \ifthenelse{\equal{\figuremode}{\fmodesvg}}{%
    \includesvg[width=\textwidth]{circle_in_air__params__nv}}{%
    \includegraphics{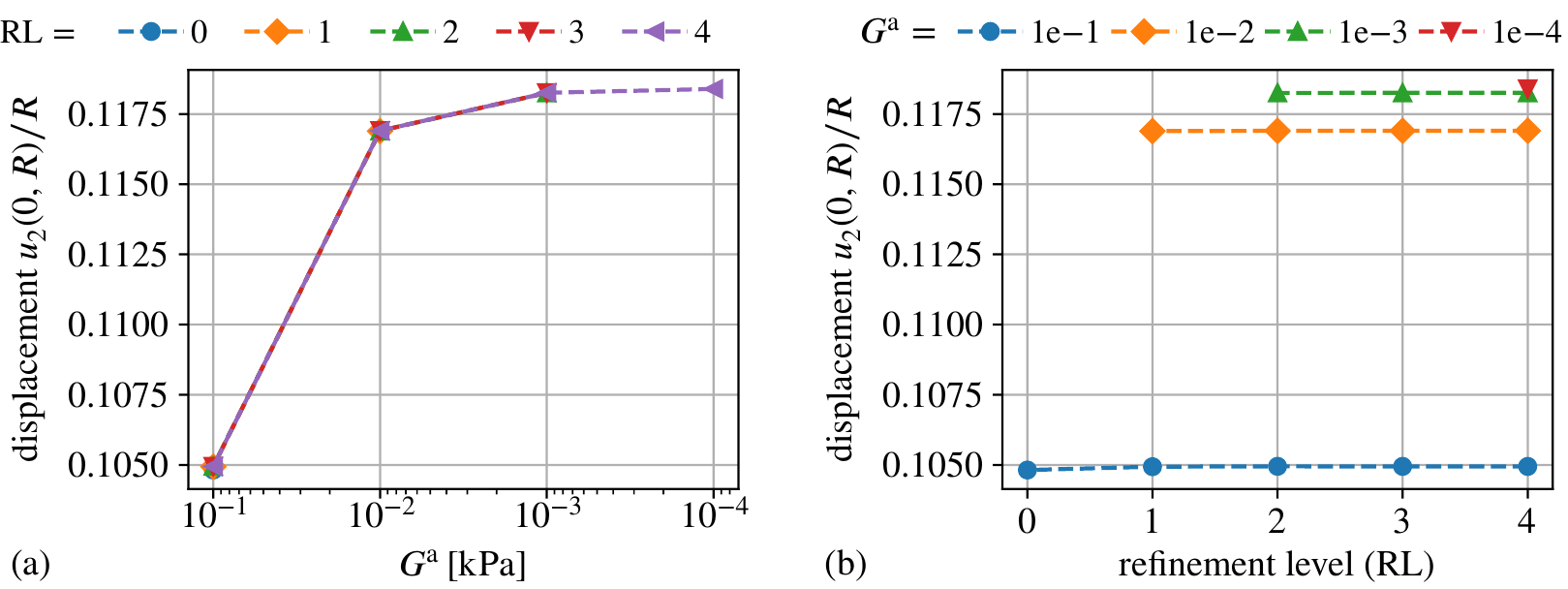}}%
    \caption{Scaled \labelu\ at $b^\infty=\SI{0.7}{\tesla}$ 
        obtained with the ``naive'' monolithic auxiliary
        stiffness approach.
        The displacement is plotted over (a)
        the auxiliary free-space stiffness $G^\text{a}$ and (b) over the mesh
        refinement level. The ``incomplete'' graphs in (a) and (b) reflect 
        the fact that the required auxiliary stiffness decreases with mesh
        refinement.}
    \label{fig:circle_in_air__params__nv}
\end{figure}

For some further insight on the underlying problem of spurious magnetic forces
we have a closer look at the capabilities of the ``naive'' fully-coupled
monolithic auxiliary stiffness approach.
Interestingly, as can be seen from Figure~\ref{fig:circle_in_air__params__nv},
the approach works quite well for very fine meshes (refinement level four (RL 4)
corresponds to more than \num{1e5} DoFs) but not so for coarse
meshes. This underlines the fact that spurious magnetic forces result from
inaccurate solutions of the magnetostatic part of the problem, which has
already been put forward in Section~\ref{sec:spurious_forces}.
Figure~\ref{fig:circle_in_air__params__nv} furthermore shows that -- if results
can be obtained at all -- convergence with respect to mesh refinement is only a
minor issue compared to the convergence with respect to the auxiliary stiffness
parameter $G^\text{a}$.
From a practical point of view, one might say ``acceptable'' results were
obtained for already for $G^\text{a}\leq\SI{1e-2}{\kilo\pascal}$ and one
refinement step. ``Accurate'' results were only obtained
$G^\text{a}\leq\SI{1e-3}{\kilo\pascal}$ requiring at least two mesh refinement
steps, i.e. going from roughly \num{1000} to \num{20000} DoFs.
A direct comparison with the methods proposed in this work is presented in
Subsection~\ref{sec:circle_in_air_convergence_all}.

\subsubsection{Maxwell traction approach via direct boundary integral computation}
\label{sec:circle_in_air_mw}

\begin{figure}[!ht]
    \centering
    \ifthenelse{\equal{\figuremode}{\fmodesvg}}{%
    \includesvg[width=\textwidth]{circle_in_air__params__mw}}{%
    \includegraphics{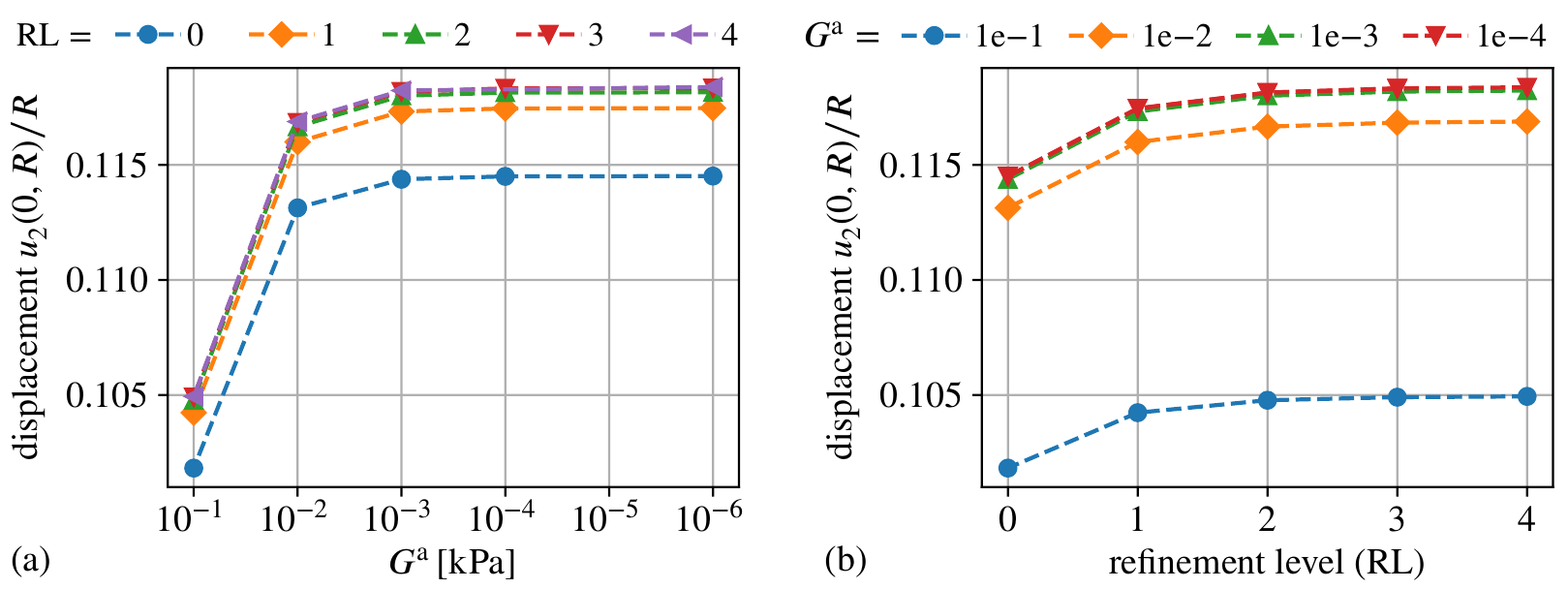}}%
    \caption{Scaled \labelu\ at $b^\infty=\SI{0.7}{\tesla}$ 
        obtained with the ``Maxwell traction'' version using boundary integrals.
        The displacement is plotted over (a)
        the auxiliary free-space stiffness $G^\text{a}$ and (b) over the mesh
        refinement level. The results for $G^\text{a}=\num{1e-4}$ and
        $\SI{1e-6}{\kilo\pascal}$ practically coincide such the graph for the
        latter is omitted in subplot (b).}
    \label{fig:circle_in_air__params__mw}
\end{figure}

Figure~\ref{fig:circle_in_air__params__mw}a
shows the effect of the auxiliary 
stiffness parameter $G^\text{a}$ for different mesh refinement levels
whereas subplot (b) shows the effect of the refinement level for given 
$G^\text{a}$. From the plots one can see that accurate results are only obtained
for $G^\text{a} \leq \SIshort{1e-3}{\kilo\pascal}$ as already observed for the
``naive monolithic'' approach. Concerning the mesh resolution, three
refinement steps (\num{86355} DoFs) yield a high-quality solution.

\subsubsection{Maxwell traction implementation via discrete force omission}
\label{sec:circle_in_air_mw_nb}

\begin{figure}[!ht]
    \centering
    \ifthenelse{\equal{\figuremode}{\fmodesvg}}{%
    \includesvg[width=0.5\textwidth]{circle_in_air__params__mw_nb}}{%
    \includegraphics{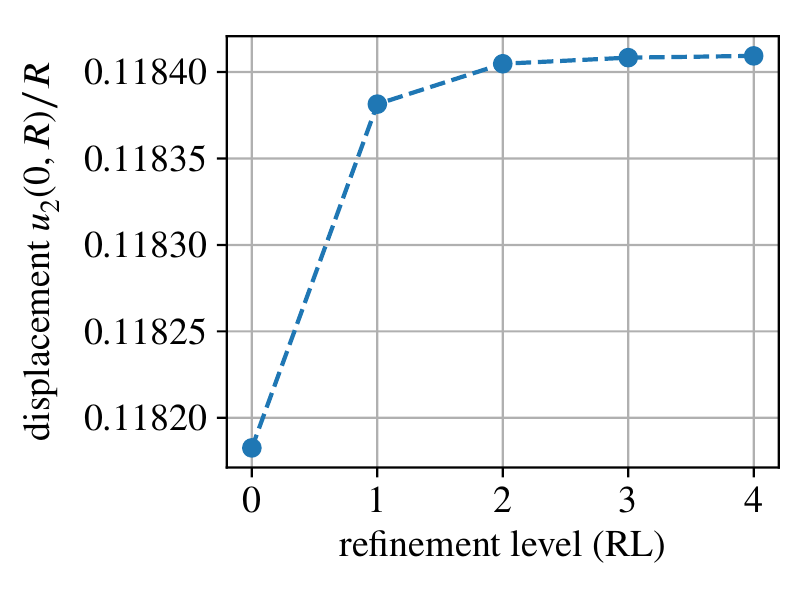}}%
        \vspace*{-5mm}
    \caption{Scaled \labelu\ at $b^\infty=\SI{0.7}{\tesla}$ for an auxiliary
        stiffness parameter $G^\text{a}=\SI{1e-6}{\kilo\pascal}$ obtained with
        the ``Maxwell traction'' version using volume integrals.
    }
    \label{fig:circle_in_air__params__mw_nb}
\end{figure}

The effect of the auxiliary stiffness parameter $G^\text{a}$ is essentially 
the same as in Figure~\ref{fig:circle_in_air__params__mw}a.
However, as shown in the mesh convergence plot depicted in
Figure~\ref{fig:circle_in_air__params__mw_nb}, 
the range of displacements is much more narrow than for the
boundary integral implementation. This indicates a much higher
accuracy of the volume-integral-based ``discrete force omission'' implementation compared
to the ``direct boundary integral'' version. 
Also, convergence with respect to mesh refinement seems to be
attained earlier with the force omission implementation.
A study of actual convergence rates is deferred to
Section~\ref{sec:circle_in_air_convergence_all}.

\subsubsection{Traction compensation implementation via direct boundary integral computation}
\label{sec:circle_in_air_tc}

\begin{figure}[!ht]
    \centering
    \ifthenelse{\equal{\figuremode}{\fmodesvg}}{%
    \includesvg[width=\textwidth]{circle_in_air__params__tc}}{%
    \includegraphics{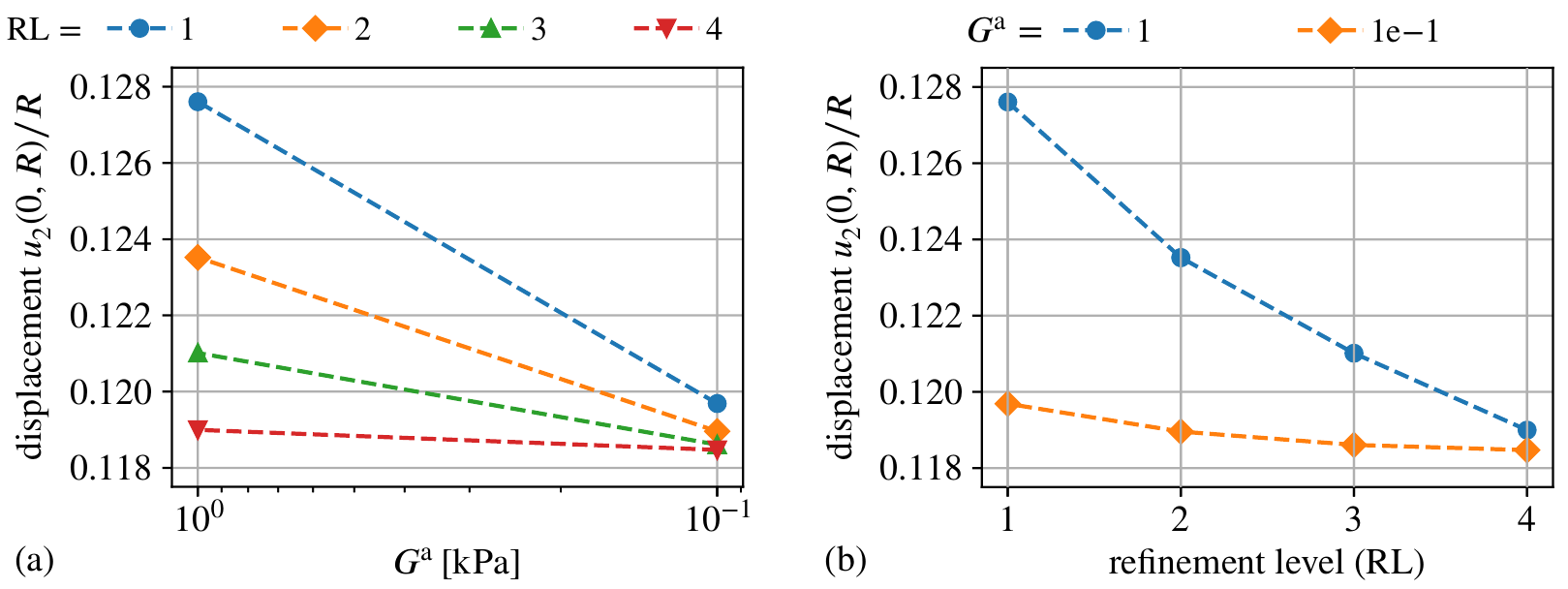}}%
    \caption{Scaled \labelu\ at $b^\infty=\SI{0.7}{\tesla}$ 
        obtained with the ``Traction compensation'' approach implemented via
        direct boundary integral computation. Subplot (a) depicts the displacement plotter over
        the auxiliary stiffness parameter $G^\text{a}$ whereas subplot (b)
        shows the displacement plotter over the mesh refinement level.
        Note that the parameter $G^\text{a}$ plays different roles in the
        ``Traction compensation'' and ``Maxwell traction'' such that a
        comparison of its effect cannot be made directly.%
    }
    \label{fig:circle_in_air__params__tc}
\end{figure}

In this case the value of the auxiliary stiffness parameter $G^\text{a}$ is
responsible for
preventing excessive deformation due to spurious magnetic forces. Therefore,
we employ considerably higher values than for the ``Maxwell traction''
implementations. Figure~\ref{fig:circle_in_air__params__tc}a 
shows a pronounced
effect of $G^\text{a}$ that clearly depends on the mesh
resolution. 
%For example, the lines for a given refinement level (mesh) become
%flatter the higher the respective refinement level (the finer the mesh) is.
Figure~\ref{fig:circle_in_air__params__tc}b indicated that the value of
the auxiliary stiffness $G^\text{a}$ is of great relevance for coarse meshes
but looses its influence for increasingly fine discretizations. 
For example, with $G^\text{a}=\SI{0.1}{\kilo\pascal}$
two refinement steps yield approximately equal results as four refinement steps
and $G^\text{a}=\SI{1.0}{\kilo\pascal}$. Thus, lower auxiliary stiffness tends
to lead to a significantly higher accuracy which is somewhat unexpected from the
theory.

\subsubsection{Traction compensation implementation via discrete force omission}
\label{sec:circle_in_air_tc_nb}

\begin{figure}[!ht]
    \centering
    \ifthenelse{\equal{\figuremode}{\fmodesvg}}{%
    \includesvg[width=\textwidth]{circle_in_air__params__tc_nb}}{%
    \includegraphics{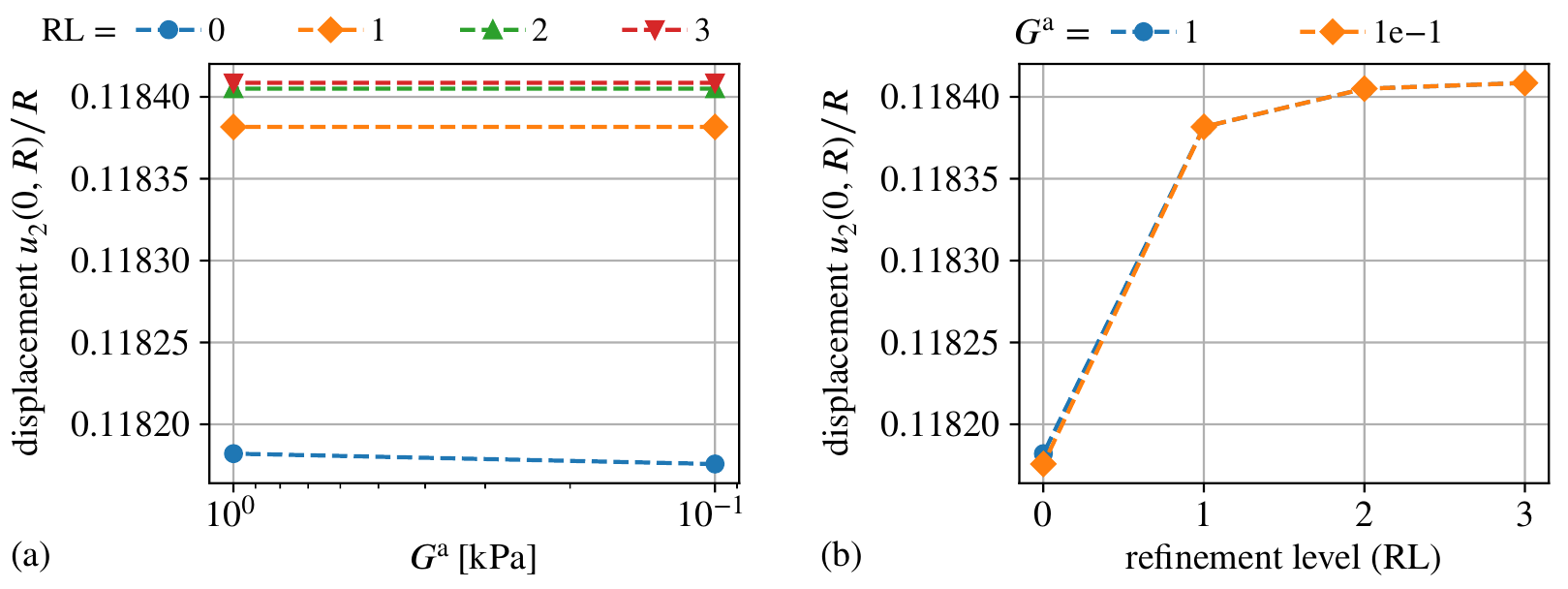}}%
    \caption{Scaled \labelu\ at $b^\infty=\SI{0.7}{\tesla}$ 
        obtained with the ``Traction compensation'' approach implemented via
        volume integrals. Subplot (a) depicts the displacement plotter over
        the auxiliary stiffness parameter $G^\text{a}$ whereas subplot (b)
        shows the displacement plotter over the mesh refinement level.}
    \label{fig:circle_in_air__params__tc_nb}
\end{figure}

As can be seen from Figure~\ref{fig:circle_in_air__params__tc_nb}a,
the force omission implementation of the traction compensation approach 
shows almost no sensitivity with respect to the auxiliary stiffness parameter
$G^\text{a}$ as expected theoretically. 
Only for the coarsest mesh (refinement level zero), on
can see a small effect in Figure~\ref{fig:circle_in_air__params__tc_nb}b.
This is due to finite deformations caused by spurious magnetic forces in the
case of $G^\text{a} = \SI{0.1}{\kilo\pascal}$. For finer meshes, the spurious
forces decline such that the solutions for 
$G^\text{a} = \SI{0.1}{\kilo\pascal}$ and
$G^\text{a} = \SI{1}{\kilo\pascal}$ practically coincide.
%This is line with the expectation that the ``Traction compensation'' approach
%converges with \emph{increasing} auxiliary stiffness. In fact, 
%in the the latter case, the auxiliary stiffness ``holds down'' the effects of
%spurious forces, which themselves decrease with increasing mesh resolution.

Also, from the tick-labels of the $y$-axis in both plots one can see that the
volume integral implementation overall is much more accurate than the
implementation based direct boundary integral computation.

\subsubsection{Comparative convergence study}
\label{sec:circle_in_air_convergence_all}

For the convergence study below, we consider as reference the solution obtained
with a staggered scheme along the lines of \citet{pelteret+etal2016}. 
In order to enhance the accuracy of their scheme, 
we keep alternating between the
coupled and the free-space adaption sub-steps until the overall solution is
converged. To avoid visual clutter, we only consider one representative set of
parameters for each method.

\begin{figure}[!ht]
    \centering
    \ifthenelse{\equal{\figuremode}{\fmodesvg}}{%
    \includesvg[width=\textwidth,pretex=%
    \newcommand{\klone}{\scriptsize{1}}%
    \newcommand{\kltwo}{\scriptsize{2}}%
    ]{circle_in_air__convergence__all}}{%
    \includegraphics{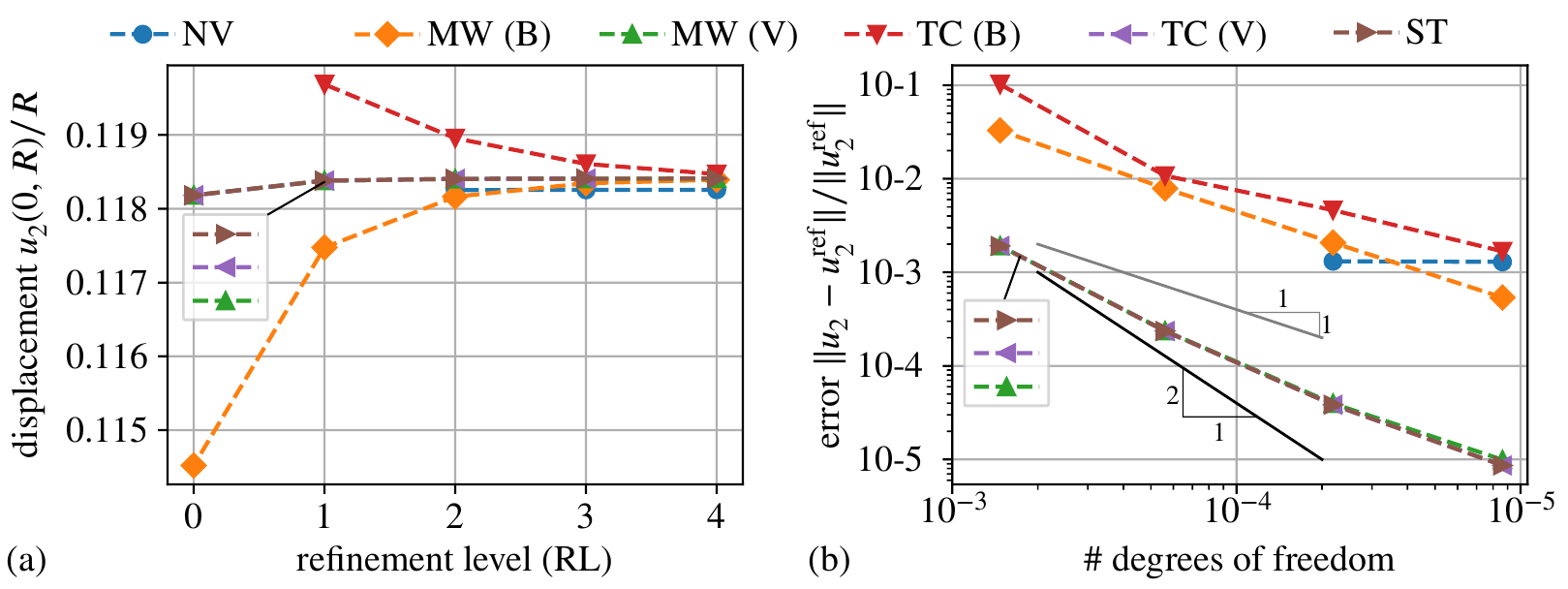}}%
    \caption{Comparative convergence study in terms of (a) the scaled
        \labelu\ plotted over the refinement level and (b) the relative
        displacement error over the number of DoFs in logarithmic
        axes. The data points in (b) correspond to refinement levels zero to
        three. The reference solution for (b) is that obtained with a staggered
        scheme for mesh refinement level four, i.e. \num{343203} degrees of
        freedom.
        The shorthand legend labels expand to:
        ``NV'' -- ``naive'' with $G^\text{a}=\SIshort{1e-3}{\kilo\pascal}$, 
        ``MW (B|V)'' -- ``Maxwell traction'' via direct \emph{B}oundary
        integrals or \emph{V}olume integrals and force omission with
        $G^\text{a}=\SIshort{1e-6}{\kilo\pascal}$, 
        ``TC (B|V)'' -- ``Traction compensation''
        via direct \emph{B}oundary integral or \emph{V}olume integral and force
        omission with $G^\text{a}=\SI{1}{\kilo\pascal}$ and 
        ``ST'' -- staggered scheme.
    }
    \label{fig:circle_in_air__convergence__all}
\end{figure}

The results for the displacement are shown in
Figure~\ref{fig:circle_in_air__convergence__all}a.
The plot confirms that the
volume integral implementations are much more accurate than their boundary
integral counterparts. Indeed, the former practically coincide with the solution
obtained with the staggered scheme.
The plot in Figure~\ref{fig:circle_in_air__convergence__all}b
depicts the convergence of the methods under consideration in terms of the
deviation from the solution of the staggered scheme at refinement level four
(\num{343203} DoFs) in dependence of the number of DoFs.
The solid gray line corresponds to linear convergence and the solid black line
to quadratic convergence in mesh resolution.
One observes that the volume integral implementations are much more accurate and
furthermore converge at super-linear rate. In contrast, the boundary integral
implementations converge at most at a linear rate.
Different from all other graphs, the one for the ``naive'' approach does not
exhibit any significant convergence because its solution (or error, respectively) is
dominated by the perturbation of the actual problem through the auxiliary
stiffness.
As a small, but nevertheless noteworthy detail we point out that only for very
fine meshes, one can see a small deviation of the ``Maxwell traction'' 
implementation via force omission (``MW (B)'') from the staggered and the ``traction
compensation'' via force omission (``TC (B)''). This stems from the auxiliary
stiffness $G^\text{a}=\SIshort{1e-6}{\kilo\pascal}$ that perturbs the
solution in the case of the ``Maxwell traction'' approach.

Due to the superior performance of the force omission implementations over their
boundary integral counterparts, the latter will be discarded from further
considerations.

\subsection{Eliminating spurious coupling in non-magnetic solid domains}
\label{sec:remedies_for_solids}

In the case of non-magnetic \emph{solid} domains that are sufficiently soft to
suffer from spurious magnetic forces, neither the ``staggered''
\citep{pelteret+etal2016} nor the non-local constraint \cite{psarra+etal2019}
approach can be applied, because both of the methods rely on the absence of
physical stiffness. In contrast, both the ``Maxwell traction'' and ``Traction
compensation'' approach can be extended to non-magnetic solids.

The adaption of the ``Maxwell traction'' approach and its implementations 
to this scenario is straight forward. 
The auxiliary strain energy density from before is simply substituted
by the actual mechanical material model.
For the ``traction compensation'' approach, the modification for solid domains
is more delicate and requires a generalization of the fundamental idea:
increase mechanical stiffness and compensate this increase by an appropriate
increase of loads.
Consider the weak balance of momentum for prescribed tractions and body force
densities
\begin{align}
    \label{eq:generic_weak_balance_of_momentum}
\int_{\volb_t} 
-\delta\defomap \cdot \left(\div\Bsigma + \bodyforce\right) \dv
+ 
\int_{\bounb_t} 
\delta{\defomap} \cdot \left(\Bsigma \cdot \bn - \bt \right) \da = 0,
\end{align}
where the (total) Cauchy stress $\Bsigma$ consists of the Maxwell stress, the
mechanical stress resulting from the actual mechanical material model.
In order to add a \emph{mechanical} zero holding down spurious
magneto-mechanical interactions, we need to be more careful than before in
order to not change the solution of the original problem.

First recall that, within a non-magnetic material, the presence of magnetic field
does not affect the mechanics. Thus, we may factor out anything related to the
Maxwell stress, i.e.
\begin{align}
\label{eq:generic_weak_balance_of_momentum_maxwell_sep_a}
&\int_{\volb_t} 
-\delta\defomap \cdot \left(\div\Bsigma^{\text{mech}} + \bodyforce\right) \dv
+ 
\int_{\bounb_t} 
\delta{\defomap} \cdot \left(\Bsigma^{\text{mech}} \cdot \bn - \bt \right) \da = 0
\sitext{and}
\label{eq:generic_weak_balance_of_momentum_maxwell_sep_b}
&\int_{\volb_t} 
-\delta\defomap \cdot \left(\div\Bsigmamw\right) \dv
+ 
\int_{\bounb_t} 
\delta{\defomap} \cdot \left(\Bsigmamw \cdot \bn \right) \da
= 0,
\end{align}
whereby \eqref{eq:generic_weak_balance_of_momentum_maxwell_sep_b} is
automatically fulfilled for any solution of the magnetostatic problem,
since the medium under consideration does not experience any magnetic forces.
Next, as the purely mechanical part
\eqref{eq:generic_weak_balance_of_momentum_maxwell_sep_a} can be treated
separately, we may scale it, i.e. multiply by some constant factor $(1 + c)$,
without affecting the solution to this equation.
For the purpose of a more compact notation we denote any quantity multiplied by
the compensation factor
$c$ with a bar, i.e. $\bar{\Bsigma}^{\text{mech}} = c{\Bsigma}^{\text{mech}}$.
Then,
\begin{align}
\label{eq:generic_weak_balance_of_momentum_maxwell_sep_c}
&\int_{\volb_t} 
-\delta\defomap \cdot 
\left[\div\left(\Bsigma^{\text{mech}} + \bar{\Bsigma}^{\text{mech}} \right) 
+ \left(\bodyforce + \bar{\bodyforce}\right) \right] \dv
\nonumber\\
& + \int_{\bounb_t} 
\delta{\defomap} \cdot \left[
    \left(\Bsigma^{\text{mech}} + \bar{\Bsigma}^{\text{mech}}\right) 
    \cdot \bn - \left(\bt + \bar{\bt}\right) \right] \da = 0,
\end{align}
and any (deformation or displacement) solution to
\eqref{eq:generic_weak_balance_of_momentum_maxwell_sep_a}
is a solution to
\eqref{eq:generic_weak_balance_of_momentum_maxwell_sep_c} and \textit{vice versa}.
Doing the usual integration by parts followed by application of the divergence
theorem yields 
\begin{align}
\label{eq:scaled_weak_balance_of_momentum_a}
&\underbrace{\int_{\volb_t} 
\Bsigma^{\text{mech}} : \grad\delta\defomap 
+ \bodyforce \cdot \delta\defomap \dv - 
\int_{\bounb_t} 
\delta{\defomap} \cdot \bt \da}_{=0}
+
\underbrace{\int_{\volb_t} 
    \bar{\Bsigma}^{\text{mech}} : \grad\delta\defomap 
    + \bar{\bodyforce} \cdot \delta\defomap \dv - 
    \int_{\bounb_t} 
    \delta{\defomap} \cdot \bar{\bt}\da}_{=0} .
\end{align}
The first set of terms are contained in the original boundary value problem whereas
the second set of terms shall be employed to suppress spurious interactions.
In order to proceed towards an implementation that is again based on omission of
discrete forces obtained by volume integrals, we have a close look at the second
``zero'' in \eqref{eq:scaled_weak_balance_of_momentum_a}.
The volume integrals can be implemented as they are, but instead of equating
them with the boundary integral exploiting $\bar{\Bsigma}\cdot\bn=\bar{\bt}$
we again just omit the resulting discrete forces corresponding to boundary
(interface) deformation, i.e. the additional terms to be considered are
\begin{align}
\label{eq:scaled_weak_balance_of_momentum_b}
0 &=
\int_{\volb_t} 
    \bar{\Bsigma}^{\text{mech}} : \grad\delta\defomap 
    + \bar{\bodyforce} \cdot \delta\defomap \dv - 
    \underbrace{\int_{\volb_t} 
    \bar{\Bsigma}^{\text{mech}} : \grad\delta\defomap 
    + \bar{\bodyforce} \cdot \delta\defomap \dv}_{\text{omit discrete ``interior''
        forces}}
\shortintertext{or}
0 &=\underbrace{\int_{\volb_t} 
    \bar{\Bsigma}^{\text{mech}} : \grad\delta\defomap 
    + \bar{\bodyforce} \cdot \delta\defomap \dv}_{\text{omit discrete boundary forces}}
.
\end{align}

In contrast to the traction compensation via discrete force omission for
air-like media, here we employ a multiple of the mechanical strain
energy density and a multiple of the body forces which must be scaled by the same
factor. Applied tractions, however, do not need any special treatment.

\subsection{Assessment of the ``Maxwell traction'' and the ``traction compensation approach for non-magnetic solid domains}
\label{sec:comparison_in_non-magnetic_solids}

In this subsection we build upon the trust in the force omission
implementations of the ``Maxwell traction'' and the ``traction compensation''
approach gained in Subsection~\ref{sec:comparison_in_air}. The direct boundary
integral implementations will not be considered due to their inferior accuracy.

Now, in the case of non-magnetic solids, we are not any longer in possession of
``reference'' solutions obtained with the staggered scheme. 
Instead, we compare the solutions of the proposed methods to results obtained
with the ``naive'' monolithic scheme, which admits direct extension to
non-magnetic solids and is quite accurate and robust as long as the 
non-magnetic solid is sufficiently stiff. In fact, what we need to demonstrate
in this section is not the accuracy of the new methods for very soft non-magnetic
media as this has already been covered in the previous section.
Instead, we want to demonstrate the correctness of the extension to actual solids.
Therefore, we may choose material parameters such that a comparison with the
``naive'' scheme is indeed reasonable.

\subsubsection{A magnetic disk in a non-magnetic carrier under gravitation-type and magnetic loading}
\label{sec:circle_in_solid}

\begin{figure}[!ht]
    \centering
    \ifthenelse{\equal{\figuremode}{\fmodesvg}}{%
    \includesvg[width=\textwidth,%
    pretex=\newcommand{\labelload}{$g=\SI{9.81}{\meter\per\second\squared}$, $b^\infty=0$}%
    \renewcommand{\labelnv}{``naive''}%
    \renewcommand{\labelmwnb}{``Maxwell traction'' (vol. int.)}%
    \renewcommand{\labeltcnb}{``Traction compensation'' (vol. int.)}%
    \newcommand{\klone}{\scriptsize{1}}%
    \newcommand{\kltwo}{\scriptsize{2}}%
    ]{circle_in_solid_disp_gravity}}{%
    \includegraphics{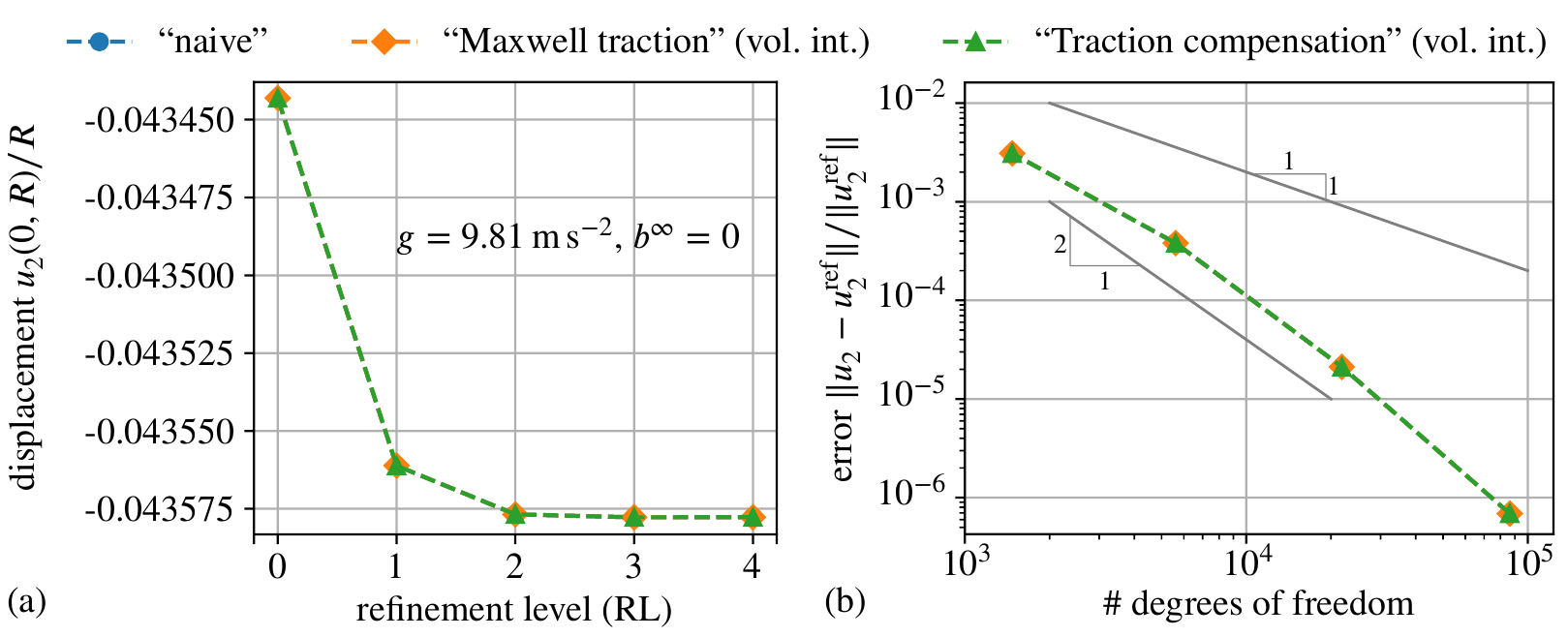}}%
    \caption{Comparative convergence study for gravitational load ($g=\SI{9.81}{\meter\per\second\squared}$, $b^\infty=0$) in terms
        of (a) the scaled
        \labelu\ plotted over the refinement level and (b) the relative
        displacement error over the number of DoFs in
        logarithmic axes. The data points in (b) correspond to refinement
        levels zero to
        three. The reference solution for (b) is that obtained with the
        ``Traction compensation'' scheme at mesh refinement level four,
        i.e. \num{343203} DoFs.
%        The shorthand legend labels expand to:
%        ``NV'' -- ``naive'' with $G^\text{a}=\SIshort{1e-3}{\kilo\pascal}$, 
%        ``MW (V)'' -- ``Maxwell traction'' via \emph{V}olume integrals with $G^\text{a}=\SIshort{1e-6}{\kilo\pascal}$ and
%        ``TC (V)'' -- ``Traction compensation''
%        via \emph{V}olume integrals with $G^\text{a}=\SI{1}{\kilo\pascal}$.
    }
    \label{fig:circle_in_solid_disp_gravity}
\end{figure}

\begin{figure}[!ht]
    \centering
    \ifthenelse{\equal{\figuremode}{\fmodesvg}}{%
    \includesvg[width=\textwidth,%
    pretex=\newcommand{\labelload}{$g=\SI{9.81}{\meter\per\second\squared}$, $b^\infty=\SI{1}{\tesla}$}%
    \renewcommand{\labelnv}{``naive''}%
    \renewcommand{\labelmwnb}{``Maxwell traction'' (vol. int.)}%
    \renewcommand{\labeltcnb}{``Traction compensation'' (vol. int.)}%
    \newcommand{\klone}{\scriptsize{1}}%
    \newcommand{\kltwo}{\scriptsize{2}}%
    ]{circle_in_solid_disp_gravity_and_field}}{%
    \includegraphics{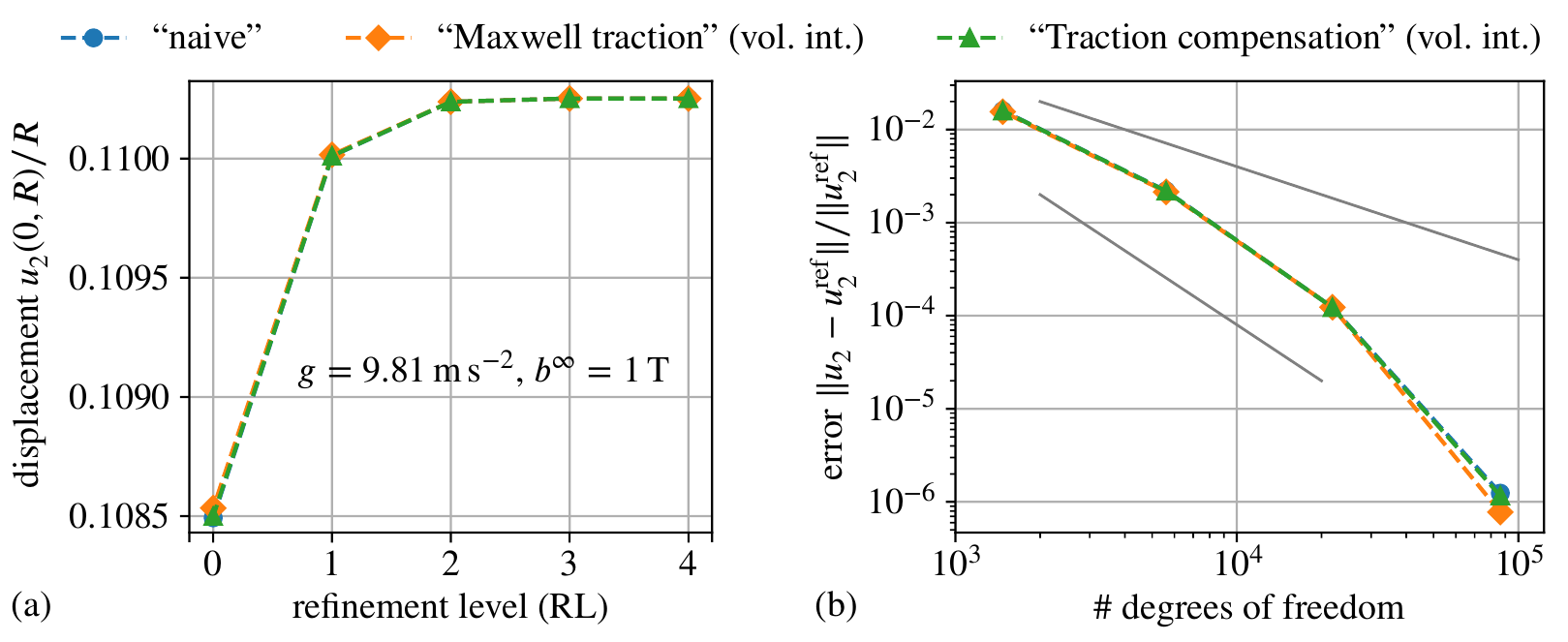}}%
    \caption{Comparative convergence study for combined gravitational and
        magnetic load ($g=\SI{9.81}{\meter\per\second\squared}$,
        $b^\infty=\SI{1}{\tesla}$) in terms of (a) the scaled
        \labelu\ plotted over the refinement level and (b) the relative
        displacement error over the number of DoFs in
        logarithmic axes. The data points in (b) correspond to refinement
        levels zero to
        three. The reference solution for (b) is that obtained with the
        ``Traction compensation'' scheme at mesh refinement level four,
        i.e. \num{343203} DoFs.
%        The shorthand legend labels expand to:
%        ``NV'' -- ``naive'' with $G^\text{a}=\SIshort{1e-3}{\kilo\pascal}$, 
%        ``MW (V)'' -- ``Maxwell traction'' via \emph{V}olume integrals with $G^\text{a}=\SIshort{1e-6}{\kilo\pascal}$ and
%        ``TC (V)'' -- ``Traction compensation''
%        via \emph{V}olume integrals with $G^\text{a}=\SI{1}{\kilo\pascal}$.
    }
    \label{fig:circle_in_solid_disp_gravity_and_field}
\end{figure}

\begin{figure}
    \centering
    \ifthenelse{\equal{\figuremode}{\fmodesvg}}{%
    \includesvg[width=\textwidth,%
    pretex=%
    \small%
    \newcommand{\lA}{\Large\textcolor{white}{A}}%
    \newcommand{\lmagn}{\Large\textcolor{red}{$\volb^\text{magn}_t$}}%
    \newcommand{\lbmagn}{\Large\textcolor{red}{$\bounb^\text{magn}_t$}}%
    \newcommand{\lbmagnref}{\Large\textcolor{red}{$\bounb^\text{magn}_0$}}%
    \newcommand{\labela}{\normalsize(a) deformed body at $g=\SI{1}{\meter\per\second\square}$, $b^\infty=0$}%
    \newcommand{\labelb}{\normalsize(b) deformed body at $g=\SI{1}{\meter\per\second\square}$, $b^\infty=\SI{1.0}{\tesla}$}%
    ]{circle_in_solid_contours}}{%
    \includegraphics{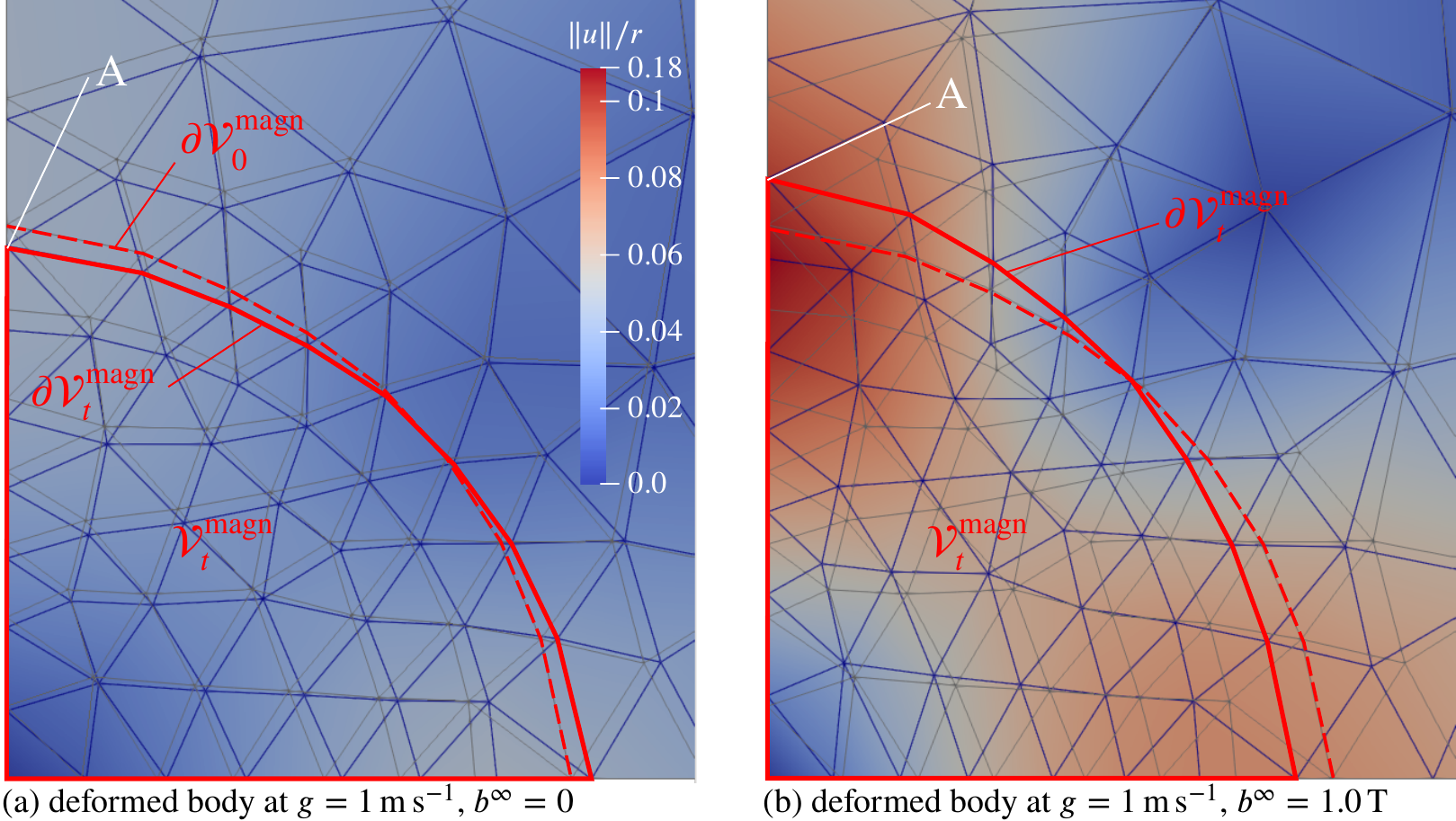}}%
    \caption{Deformed configurations under (a) purely gravitational load and (b) combined gravitational and
        magnetic loads. The solid red line indicates the deformed magnetic domain
        $\volb^\text{magn}_t$ whereas the dashed red line indicates the boundary of the undeformed magnetic
        domain $\bounb^\text{magn}_0$. The contours refer to the displacement magnitude, whereby the colors
        are adjusted to displacement range in (b). 
        Please note that while the mesh nodes are connected by straight lines in the figure, simulations
        actually have been carried out with second-order geometry descriptions.}
    \label{fig:circle_in_solid_contours}
\end{figure}

Here we in essence reuse the boundary value problem from
Sections~\ref{sec:spurious_forces} and \ref{sec:comparison_in_air} but
with a specific energy density assigned to the non-magnetic domain.
For simplicity, we again choose an energy density of the form
\eqref{eq:demo_magnetic_material} with 
parameters $\{G, G', \chi\} = \{\SI{1}{\kPa}, \SI{50}\kPa, 10\}$ for the
magnetic domain and
$\{G, G', \chi\} = \{\SI{0.5}{\kPa}, \SI{25}\kPa, 0\}$ for the
non-magnetic domain.
In addition, we assign mass densities\footnote{Under assumption of unit length ``thickness''.} to both the magnetic and the non-magnetic
solid leading to gravitation-like forces. To be specific, the mass density in the
magnetic domain is $\rho_0^\text{magn} = \SIshort{1000}{\kg\per\meter\cubed}$
whereas in the non-magnetic domain we employ
$\rho_0^\text{nonm} = \SIshort{100}{\kg\per\meter\cubed}$. Gravity points
in negative vertical ($y$) direction such that the resulting body forces are given as
\begin{align}
    \label{eq:body_forces}
    &\bodyforce{}^\text{|magn}_0 = \left(0, -g\,\rho^\text{magn}_0\right)^\text{T}
    &\text{and}&
    &\bodyforce{}^\text{|nonm}_0 = \left(0, -g\,\rho^\text{nonm}_0\right)^\text{T}
\end{align}
where $g$ denotes the gravitational loading
parameter.%
    \footnote{The example serves for the evaluation of the capabilities of the
        methods under consideration. Thus, we allow ourselves to deviate from
        actual physical parameters.}
In Figure~\ref{fig:circle_in_solid_disp_gravity}a 
we show the scaled vertical
displacement $u_2(0, R)/R$ at $g=\SIshort{1}{\meter\per\second\square}$ and
$b^\infty=0$, i.e. purely gravitational loading. 
The corresponding mesh convergence plot is shown in
Figure~\ref{fig:circle_in_solid_disp_gravity}b.
In both subplots one can observe the perfect agreement of all three
methods.
The picture is similar in
Figure~\ref{fig:circle_in_solid_disp_gravity_and_field},
where the scaled vertical
displacement $u_2(0, R)/R$ is depicted at the magneto-mechanical loading
state
$g=\SIshort{1}{\meter\per\second\square}$ and
$b^\infty=\SIshort{1}{\tesla}$.
The deformed magnetic bodies at purely gravitational and combined loading
are depicted in subplots (a) and (b), respectively, of
Figure~\ref{fig:circle_in_solid_contours}.
There one can see that displacement of point ``A'' is downwards under gravitational load (subplot a) and
upwards under combined loading (subplot b).

\subsubsection{A bilayer beam under gravity air magnetic field embedded in an air-like domain}
\label{sec:bilayer_beam}

\begin{figure}
    \centering
    \ifthenelse{\equal{\figuremode}{\fmodesvg}}{%
    \includesvg[width=0.9\textwidth,%
                pretex=%
                \newcommand{\lbext}{$\bb^\infty$}%
                \newcommand{\lgravity}{$g$}%
                \newcommand{\labelL}{$L$}%
                \newcommand{\labelR}{$l=L/10$}%
                \newcommand{\labelr}{$l/10$}%
                \newcommand{\labelrhalf}{$l/20$}%
                \newcommand{\labelLhalf}{$L/2$}%
                \newcommand{\lairdomain}{$\volfree$}%
                \newcommand{\magndomain}{$\volb^\text{magn}$}%
                \newcommand{\nonmagndomain}{$\volb^\text{nonm}$}%
                ]{bilayer_beam_example}}{%
    \includegraphics{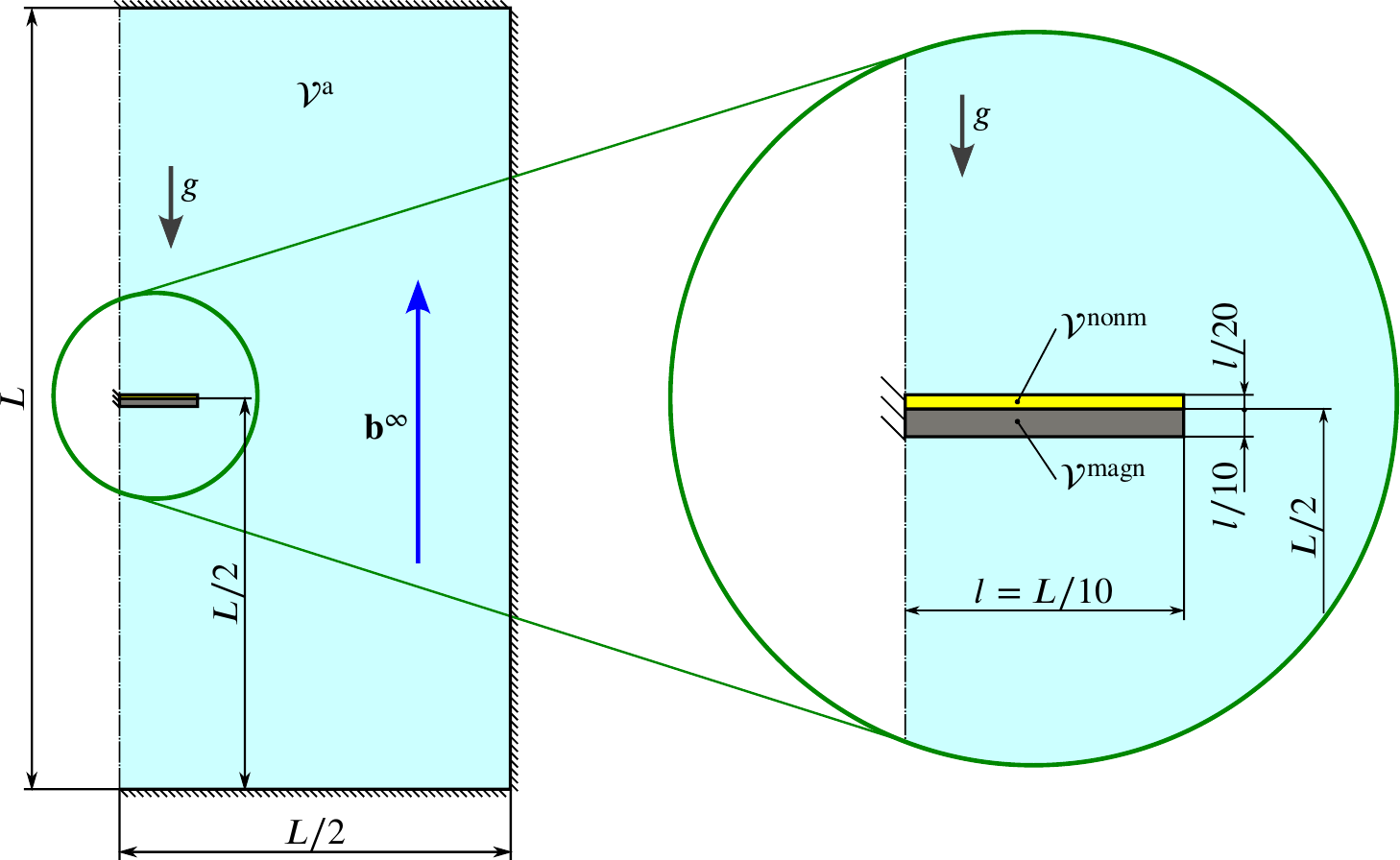}}%
    \caption{Boundary value problem of a beam that consists of a magnetic $\volb^\text{magn}$
             and a non-magnetic layer $\volb^\text{nonm}$. It is clamped at its vertical symmetry axis and
             surrounded by ``empty'' space $\volfree$. The bilayer beam is exposed to gravity
             $g$ and a uniform external magnetic field $\bb^\infty$.}
    \label{fig:bilayer_beam_example}
\end{figure}

As final example that demonstrates the capabilities of the proposed schemes
we consider a beam consisting of a magnetic and non-magnetic layer that is embedded
in an air-like domain as depicted in Figure~\ref{fig:bilayer_beam_example}.
It depicts a bilayer beam consisting of a magnetic $\volb^\text{magn}$
and a non-mangetic layer $\volb^\text{nonm}$ that is clamped at its vertical symmetry axis and
surrounded by ``empty'' space $\volfree$. This setting is inspired by prior examples of MRE beams 
\citep{zhao+etal2019,mukherjee+rambausek+danas2021,rambausek+mukherjee+danas2022,moreno+etal2022}
and the prospective application of magnetoactive materials as mechanically active substrate for 
biological experiments \citep{gonzalez-rico+etal2021}.
In the simulations, the beam is first exposed to gravity
$g=\SIshort{9.81}{\meter\per\second\squared}$ 
and in a second step to a combined loading through gravity and a uniform external magnetic field
$\bb^\infty=\SI{1}{\tesla}$.
The generic energy density function employed is given as
\begin{align}
\PsiCH &= \frac{G}{2} \left[\Tr\bC - 2\log J - 3\right] + 
\frac{G'}{2} \left(J - 1\right)^2
\nonumber \\
&\quad
-\frac{J \mu_0 (\msat)^2}{\chi} \log\left[\cosh\left(\frac{\chi\sqrt{\bH \cdot (\bC^{-1} \cdot \bH)}}{\mu_0\msat}\right)\right]
+ \PsiCHvac .
\end{align}
It has the same purely mechanical neo-Hookean type contribution as in previous
examples. The magnetic part, however, now models saturation
$\norm{ \bbm} \to \msat$ for $\norm{\bh}\to\infty$.

The material parameters for this example are collected in
Table~\ref{tab:bilayer_beam_example_params}.
\begin{table}[!ht]
    \centering
    \caption{Material parameters for the solid domains in the bilayer beam example}
    \label{tab:bilayer_beam_example_params}
    \begin{tabular}{l c c c c c }
        \toprule[0.35mm]
        Domain & $G\perunit{\kilo\pascal}$ & $G'\perunit{\kilo\pascal}$ & 
        $\rho\perunit{\kilogram\per\meter\cubed}$ & $\chi\perunit{\unitone}$ 
        & $\msat\perunit{\mega\ampere\per\meter}$ \\
        \midrule[0.25mm]
        $\volb^\text{magn}$ & \num{2e3} & \num{100e3} & \num{2e3} & 10 & 1 \\
        $\volb^\text{nonm}$ & \num{1} & \num{50} & \num{1e3} & 0 & -- \\
%        \midrule[0.25mm]
%        (compensation) & $G^\text{a}\perunit{\kilo\pascal}$ & $G^\text{a}'\perunit{\kilo\pascal}$ & & & \\
%        \midrule[0.25mm]
%        $\volfree$ (MW) & \num{1e-6} & \num{5e-5} & 0 & 0 & -- \\
%        $\volfree$ (TC) & \num{1e2} & \num{5e-3} & 0 & 0 & -- \\
        \bottomrule[0.35mm]
    \end{tabular}
\end{table}
They have been chosen with great care to render a physically reasonable 
example that challenges the numerical methods under evaluation.
%The shorthands ``MW'' and ``TC'' stand for ``\emph{M}axwell \emph{T}raction''
%and ``\emph{T}raction \emph{C}ompensation''.
Furthermore, when the ``Maxwell
traction'' approach is applied, the air-like domain $\volfree$ is equipped with
$G^\text{a}=\SIshort{1e-6}{\kilo\pascal}$. 
In the case of ``traction
compensation'', both the non-magnetic and air-like domains have
$G^\text{a}=\SIshort{1e2}{\kilo\pascal}$, which amounts to a compensation factor of
$c=100$ in the solid domain $\volb^\text{nonm}$. This means that the 
auxiliary stiffness added to the solid is much higher that the actual stiffness.
This deliberate choice helps us to confirm that all effects from the additional
stiffness are compensated properly.

\begin{figure}[!ht]
    \ifthenelse{\equal{\figuremode}{\fmodesvg}}{%
    \includesvg[width=\textwidth,
    pretex=\small%
    \newcommand{\lmagn}{\textcolor{red}{$\volb^\text{magn}$}}%
    \newcommand{\lnonmagn}{\textcolor{red}{$\volb^\text{nonm}$}}%
    \newcommand{\labela}{\normalsize(a) deformed body at $g=\SI{9.81}{\meter\per\second\square}$, $b^\infty=0$}%
    \newcommand{\labelb}{\normalsize(b) deformed body at $g=\SI{9.81}{\meter\per\second\square}$, $b^\infty=\SI{1.0}{\tesla}$}%
    ]{bilayer_beam_contours}}{%
    \includegraphics{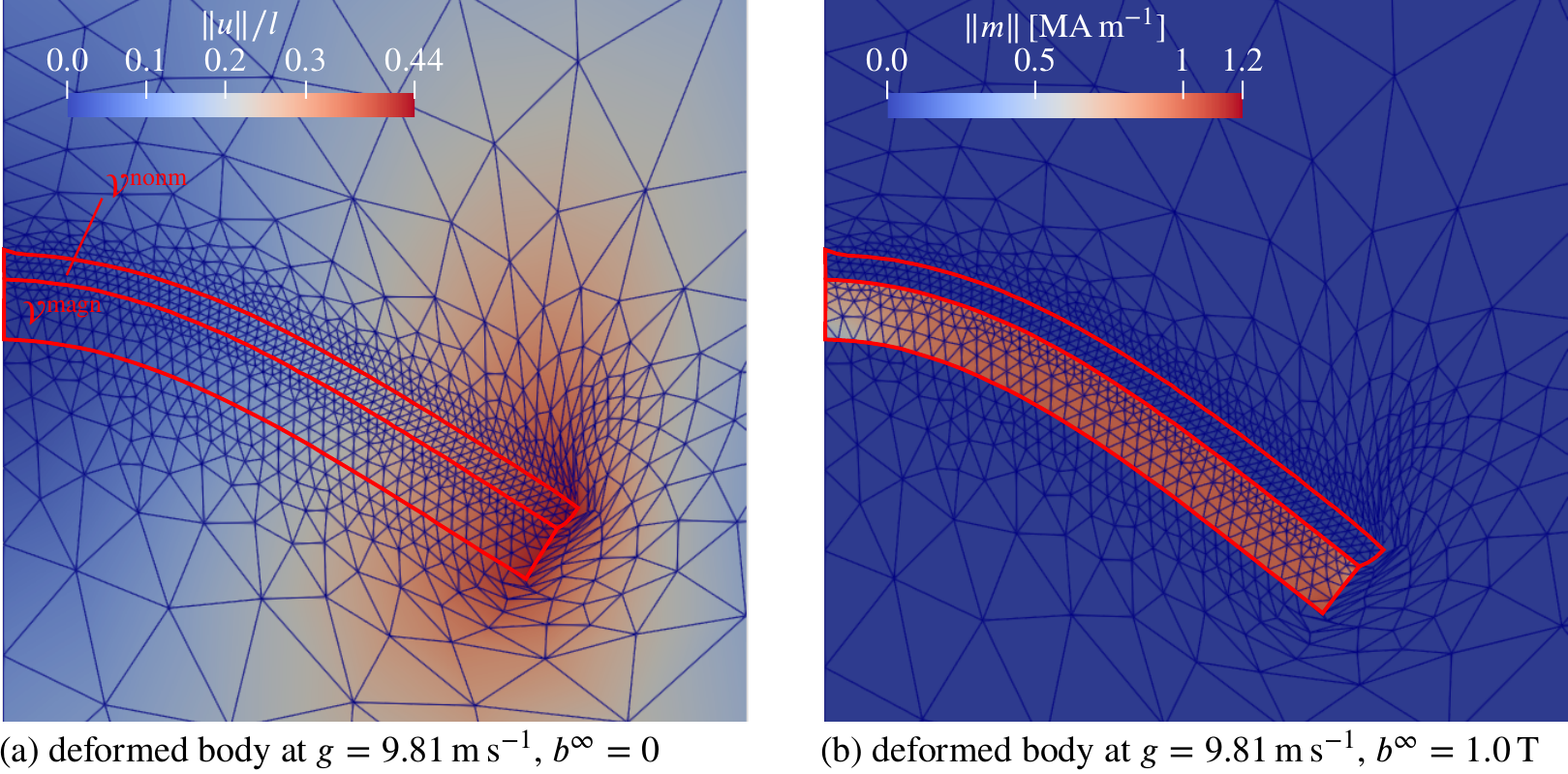}}%
    \label{fig:bilayer_beam_contours}
\end{figure}

In contrast to all foregoing examples, the geometry of the bilayer beam has sharp
corners. This complicates convergence studies because the ``naive'' scheme has
severe problems in such cases. The reason are magnetic field concentrations
near the corners leading to pronounced spurious magnetic forces that cannot be
simply alleviated with mesh refinement
\citep[Chapter~9, Section~3.1.3]{rambausek2020}. 
Therefore, we omit a rigorous convergence study
and instead show contours of solution states in
Figure~\ref{fig:bilayer_beam_contours}.
In subplot (a) we observe a gravity-driven bending-dominated deformation of the magnetic layer
$\volb^\text{magn}$ and a slightly more general deformation pattern of the non-magnetic layer
$\volb^\text{nonm}$. The color contours in (a) correspond to the deformation magnitude.
Subplot (b) shows the deformed configurations under combined loading. The deformation
modes are very similar to those under purely gravitational loading, the displacements
magnitudes are further increased through the applied magnetic field. The color contours in
(b) correspond to the magnetization magnitude. They first of all confirm that only the
magnetic layer exhibits magnetization as expected and, second, that the magnetization
distribution is far from being perfectly uniform or trivial. This underlines the importance
of full-field simulations of such boundary value problems.
We highlight that in both subplots of Figure~\ref{fig:bilayer_beam_contours} one can
observe severe deformations of the ``empty'' domain surrounding the beam in the vicinity of
the tip of the beam. While this was not of a great concerns here, in practice one could
``smoothen'' the deformation field by carefully increasing the auxiliary stiffness in such regions
relative to that in other areas of the ``empty'' domain. Please also note that the
post-processed deformation is only of first order but the actual deformation is the
simulation is of second order, such that severely distorted elements in the vicinity of the
beam tip are not displayed exactly how they actually appear in the simulations.
In any case, the solutions obtained for the Maxwell traction and the traction
compensation scheme coincide almost perfectly.
In this context we highlight that we did not need to fine-tune the auxiliary
stiffness parameters to achieve excellent agreement between both approaches.
This underlines not only their computational robustness but also their
robustness with respect to the methods parameters.
In general, we for the Maxwell traction scheme recommend to choose the
auxiliary stiffness parameter in a range of $10^{-6}$ of the softest solid
parameter.
For the ``Traction compensation'' we recommend stiffness parameters roughly in
the range of those employed for the magnetic solids.

\section{Conclusion}
\label{sec:conclusion}

This work is centered around the issue of spurious magneto-mechanical
interactions that is pervasive in ``naive'' fully-coupled numerical simulations
of deformable magnetic bodies such as MREs.
The key contributions are, first, a thorough characterization of
the underlying issue. Second, we present two novel approaches that effectively
eliminate or suppress spurious
coupling in both vacuum- or air-like media and non-magnetic solid. 
The first scheme relies on the so-called ``Maxwell tractions'' that removes all
unwanted magneto-mechanical interaction from the interior of the non-magnetic domain.
For definiteness of the solution in ``empty'' (air-like or 
vacuum domains) the method, only needs a negligibly small auxiliary stiffness even
in comparison with very soft solids.
The second scheme is somewhat dual to the first in that it
employs a sufficiently large auxiliary mechanical stiffness to suppress unphysical
magneto-mechanically driven deformation in non-magnetic domains.
The additional stiffness is balanced by additional body force contributions
and a removal or a compensation of the resulting tractions on the interface 
to neighboring bodies. The common advantages of our proposed methods over
existing successful schemes, in particular in comparison with the staggered
schemes \citep{pelteret+etal2016,liu+etal2020} and monolithic schemes based on
non-local constraints \citep{psarra+etal2019,rambausek+mukherjee+danas2022} are
twofold. The first is the ease of implementation atop of ``naive'' monolithic
FE simulations because, as has been shown, the successful volume-integral-based 
versions rely on only slight modifications of the weak form. There is no need 
for modifications to the overall solution
procedure, nor cumbersome preprocessing and modifications to the sparsity
structure of the linearized system matrix.
Second, as demonstrated successfully, both of our approaches are not only applicable
for air-like environments but also for actual, very soft non-magnetic \emph{solids}.
Another advantage over the staggered scheme is that both of the proposed
approaches directly allow for the consistent linearization of the fully coupled
system, which positively affects the convergence of nonlinear solvers.
A mild disadvantage is that the resulting linear systems are non-symmetric, which
increases the effort required for linear solves. To our experience, this is
usually outweighed by the improved convergence and robustness of the methods.

We have assessed implementations of both proposed methods regarding their accuracy
and their effectiveness regarding the elimination of spurious coupling in
comparison with existing approaches. 
In particular, we have shown the convergence of the ``Maxwell traction'' approach
for decreasing auxiliary stiffness, while the ``traction compensation'' method
converges for increasing auxiliary stiffness.
In this context is important to note that while there remains a parameter to be set, 
namely the auxiliary stiffness in one form or another, there is typically no need 
for tuning this parameter with great care. For the ``Maxwell traction'' one may go right in the
range of $10^{-6}$ of the softest solid parameter, for the ``Traction compensation'' one may
simply employ some stiffness in the range of the magnetic solids under consideration.
The critical comparison with existing schemes has demonstrated the competitiveness
of the force-omission-based variants of the approaches proposed, with
small advantages for the ``traction compensation'' approach.
Thus, in view of the ease of computer implementation and simple choice of parameters
both of the proposed methods can hope for wide adoption in future numerical
investigations on MREs and related materials 
or problems, e.g. in electromechanics.

\section*{Acknowledgments}

\noindent
Both of the authors acknowledge financial support from the Austrian Science Fund
(FWF) project F 65.

\noindent
For the purpose of open access, the author has applied a CC BY public copyright licence to
any Author Accepted Manuscript version arising from this submission.

\noindent
Moreover, we would like to thank \href{https://zenodo.org}{zenodo} for hosting the
\href{\suppdoi}{supplementary material}.

\section*{References}
\bibliographystyle{ascarticle-harv}
\bibliography{references}

%% APPENDIX

% \appendix

\end{document}